\documentclass[11pt]{amsart}
\usepackage[T1]{fontenc}
\usepackage{lmodern,microtype}
\usepackage[margin=1.1in]{geometry}
\usepackage{amsmath,amssymb,mathtools,enumitem}
\usepackage{tikz}
\usetikzlibrary{arrows.meta,calc,positioning,decorations.pathreplacing,patterns}
\usepackage{xurl}
\usepackage{placeins,needspace}
\usepackage[colorlinks=true,linkcolor=blue!45!black,citecolor=blue!45!black,urlcolor=blue!45!black]{hyperref}
\usepackage[capitalize,nameinlink]{cleveref}

\newtheorem{theorem}{Theorem}[section]
\newtheorem{proposition}[theorem]{Proposition}
\newtheorem{lemma}[theorem]{Lemma}
\newtheorem{corollary}[theorem]{Corollary}
\newtheorem{mainthm}{Theorem}

\theoremstyle{definition}
\newtheorem{definition}[theorem]{Definition}

\theoremstyle{remark}
\newtheorem{remark}[theorem]{Remark}
\numberwithin{equation}{section}
\crefname{mainthm}{Theorem}{Theorems}
\crefname{equation}{}{}
\crefname{figure}{Figure}{Figures}
\hypersetup{pdftitle={The Borel Ruziewicz problem},pdfauthor={Konstantin Slutsky}}

\newcommand{\R}{\mathbb R}

\newcommand{\Z}{\mathbb Z}
\newcommand{\Bor}{\mathcal B}
\newcommand{\ind}{\mathbf 1}
\newcommand{\SO}{\operatorname{SO}}
\newcommand{\TV}{\mathrm{TV}}
\newcommand{\dist}{\operatorname{dist}}
\newcommand{\Leb}{\mathcal L}
\newcommand{\Div}{D_S}
\newcommand{\Bb}{B_b}
\newcommand{\eps}{\varepsilon}

\setlist{itemsep=2pt,topsep=4pt}
\title{The Borel Ruziewicz problem}
\author{Konstantin Slutsky}
\address{Department of Mathematics, Iowa State University, 411 Morrill
Road, Ames, IA 50011, USA}
\email{kslutsky@iastate.edu}
\thanks{The results of this paper were obtained by AI; see the
paragraph \protect\hyperref[par:ai]{``Use of AI''} at the end of the introduction.}
\thanks{This material is based upon work supported by the National Science
Foundation under CAREER Grant No.~DMS-2540790.}


\begin{document}

\begin{abstract}
We prove that for $n\ge2$ the normalized surface measure is the only
rotation-invariant finitely additive probability measure on the Borel
subsets of $\mathbb S^n$, and that for $d\ge3$ Lebesgue measure is the
only finitely additive measure on the bounded Borel subsets of $\R^d$
which is invariant under orientation-preserving isometries and normalized
on the unit cube. This answers the Borel version of the Ruziewicz
problem. On $\mathbb S^n$ it suffices to assume invariance under a finite
set of rotations with a spectral gap, and on $\R^d$ under two such sets
of rotations about two distinct centres. The main step is a criterion for
invariant finitely additive probability measures to vanish on null sets.
For a finite nonempty symmetric set of homeomorphisms of a
$\sigma$-compact metrizable space with a quasi-invariant Borel measure $m$
which is finite on compact sets, such as Lebesgue measure on $\R^d$,
every
invariant finitely additive probability measure on Borel sets vanishes
on $m$-null sets if and only if there is $\delta>0$ such that every
countably additive Borel probability measure supported on an $m$-null
Borel set is moved by some generator by at least $\delta$ in total
variation. Sufficiency is proved by constructing, for every invariant
null Borel set, a bounded Borel field on the orbit graph, supported on
the set and with divergence at least one there.
\end{abstract}

\maketitle

\section{Introduction}\label{sec:intro}

The Ruziewicz problem asks whether Lebesgue measure is characterized by
finite additivity and invariance under isometries. For the unit sphere
$\mathbb S^n$ let
$\sigma_n$ be the normalized surface measure. Is
$\sigma_n$ the only finitely additive, rotation-invariant probability
measure defined on all Lebesgue measurable subsets of $\mathbb S^n$?
Banach showed that the answer is negative for the circle~\cite[Th\'eor\`eme~II]{Banach}. For $n\ge2$ the answer is positive. This
was proved for $n\ge4$ by Margulis~\cite{Margulis} and
Sullivan~\cite{Sullivan}, and for $n=2,3$ by Drinfeld~\cite{Drinfeld}. For Euclidean space
$\R^d$, $d\ge3$, Margulis proved the analogous uniqueness theorem for
finitely additive functions on the ring of bounded Lebesgue measurable
sets, invariant under all isometries and normalized on the unit cube~\cite{MargulisEucl}.
Boutonnet and Ioana~\cite{BI} extended this to
countable dense subgroups $\Gamma$ of the orientation-preserving
isometries of $\R^d$, $d\ge3$, whose image in $\SO(d)$ acts on $\SO(d)$
with a spectral gap.

For $n\ge2$ the spherical results are proved with invariant means, in
two steps. First, a rotation-invariant finitely additive probability
measure on the Lebesgue sets vanishes on null sets. This is an
observation of Tarski based on paradoxical decompositions (see~\cite[footnote~1, p.~122]{Sullivan}), and its proof uses arbitrary
subsets of null sets. It follows that such a measure is the same thing
as an invariant mean on the space $L^\infty(\sigma_n)$ of essentially
bounded functions modulo null functions, that is, a positive
normalized invariant linear functional on $L^\infty(\sigma_n)$. Second,
one finds finitely many rotations whose action on the space
$L^2(\sigma_n)$ of square-integrable functions has a
spectral gap, and shows that a spectral gap forces integration to be the
only invariant mean on $L^\infty$ (see~\cite[Theorem~1.1(a)]{LiNg} for a
precise modern statement). Sarnak's book~\cite[Chapter~2]{Sarnak}
gives an account of this solution.

In this paper we consider finitely additive measures defined only on \emph{Borel} sets. This leads
to the \emph{Borel Ruziewicz problem}~\cite[Question~13.15]{Wagon}. The Borel $\sigma$-algebra is
smaller than the Lebesgue $\sigma$-algebra, so a priori there could be more invariant finitely
additive measures on it, and the classical route does not apply directly. Tarski's argument uses
subsets of null sets which need not be Borel, and passing from a finitely additive measure on Borel
sets to a functional on $L^\infty$ requires knowing that it vanishes on Borel null sets. Grabowski,
M\'ath\'e and Pikhurko~\cite{GMP} proved, for $n\ge3$, that every
nonnegative finite-valued finitely additive function on the relatively compact Lebesgue measurable
sets with the Baire property (sets differing from an open set by a meagre set) which is invariant
under $\SO(n)$ acting on $\mathbb S^{n-1}$, or under the orientation-preserving isometries of
$\R^n$, is a scalar multiple of the standard measure. Mycielski and Tomkowicz~\cite{MT}
gave another proof of uniqueness on the ring of Lebesgue measurable sets with the Baire property in
$\mathbb S^n$ and on the ring of such bounded sets in $\R^{n+1}$, $n\ge2$, assuming invariance under
isometries. A finitely additive measure defined only on Borel sets need
not extend to these larger rings, so these results do not settle the
Borel case. Dougherty and
Foreman~\cite{DF} showed that no rotation-invariant finitely
additive probability measure on the Borel subsets of $\mathbb S^2$
vanishes on all meagre sets, and Marks and Unger~\cite{MUBaire}
showed that every paradoxical action of a group by Borel automorphisms
of a Polish space has a paradoxical decomposition whose pieces have the
Baire property. Finitely additive and countably additive invariant
measures on Borel sets can also behave quite differently. Conley,
Jackson, Marks, Seward and Tucker-Drob~\cite{CJMST}
constructed a continuous free action of $(\Z/2\Z)^{*3}$ on a Polish space
which admits an invariant finitely additive probability measure on the
Borel sets but no invariant countably additive one.

\subsection{Main results}
A \emph{probability charge} on a measurable space $(\Omega,\mathcal A)$
is a finitely additive probability on $\mathcal A$:
\[
 \mu\colon\mathcal A\to[0,1],\qquad\mu(\Omega)=1,\qquad
 \mu(A\cup B)=\mu(A)+\mu(B)\quad(A\cap B=\varnothing).
\]
Let $S$ be a finite nonempty symmetric ($S=S^{-1}$) set of
measure-preserving automorphisms of a probability space $(\Omega,m)$. We
say that $S$ has a \emph{spectral gap} (for $m$) if there is $\kappa>0$
with
\[
 \max_{s\in S}\|u\circ s^{-1}-u\|_{L^2(m)}\ge\kappa\,\|u-m(u)\|_{L^2(m)}
 \qquad\bigl(u\in L^2(m)\bigr),
\]
where $m(u)=\int u\,dm$. For rotations of $\mathbb S^n$ we take
$m=\sigma_n$. Finite sets of rotations with a spectral gap exist in every
dimension $n\ge2$ by the aforementioned results of Margulis, Sullivan and Drinfeld, and by
later work of Bourgain and Gamburd and of Benoist and de Saxc\'e they are
in fact plentiful (\cref{rem:many-gap-sets}).

Our main result is the positive solution to the Borel Ruziewicz problem.

\begin{mainthm}\label{thm:sphere}
Let $n\ge2$, and let $\sigma_n$ be the normalized surface measure on
$\mathbb S^n$. Then $\sigma_n$ is the only $\SO(n+1)$-invariant
probability charge on the Borel subsets of $\mathbb S^n$. More generally,
if $S\subseteq\SO(n+1)$ is a finite nonempty symmetric set with a
spectral gap, then $\sigma_n$ is the only $S$-invariant probability
charge on the Borel $\sigma$-algebra $\Bor(\mathbb S^n)$.
\end{mainthm}

\begin{mainthm}\label{thm:euclid}
Let $d\ge3$, and let $\mu$ be a finitely additive function from the
bounded Borel subsets of $\R^d$ to $[0,\infty]$ which is invariant under
all orientation-preserving isometries of $\R^d$ and satisfies
$\mu([0,1]^d)=1$. Then $\mu(A)$ is the Lebesgue measure of $A$ for every
bounded Borel set $A$. More generally, the same conclusion holds if
$\mu$ is only invariant under the rotations about $p$ from a finite set
$S^{(1)}\subseteq\SO(d)$ and the rotations about $q$ from a finite set
$S^{(2)}\subseteq\SO(d)$, where $p,q\in\R^d$ are distinct and $S^{(1)},S^{(2)}$ are
nonempty, symmetric and have spectral gaps on $\mathbb S^{d-1}$.
\end{mainthm}

Both theorems fail in lower dimensions. Banach constructed a
rotation-invariant finitely additive probability measure on all subsets
of the circle which gives mass one to some Lebesgue null set $W$~\cite[Th\'eor\`eme~II]{Banach}. It also gives mass one to any Borel null
set containing $W$, so its restriction to Borel sets is not $\sigma_1$.
Similarly, Banach~\cite[Th\'eor\`eme~IV*]{Banach} constructed an
isometry-invariant finitely additive function on the bounded subsets of
the plane, normalized on the unit square, with value one on a bounded
null set. So \cref{thm:euclid} fails for $d=2$.

The main new ingredient is a general criterion for invariant charges to
vanish on null sets. Let $X$ be a $\sigma$-compact metrizable space
(a countable union of compact sets), $S$ a finite nonempty symmetric set
of homeomorphisms of $X$, and $m$ a Borel measure on $X$ which is finite
on compact sets and whose null sets are preserved by $S$. We say
that $(X,S,m)$ has the \emph{singular displacement gap} if there is
$\delta>0$ such that every Borel probability measure $\nu$ supported on
an $m$-null Borel set, written $\nu\perp m$, is moved by some generator
by at least $\delta$ in total variation:
\[
 \max_{s\in S}\|s_*\nu-\nu\|_\TV\ge\delta\qquad\text{for all Borel probability measures }\nu\perp m.
\]
Here $s_*\nu(A)=\nu(s^{-1}A)$, and $\|\cdot\|_\TV$ is the total
variation norm, normalized so that mutually singular probability
measures are at distance two.

\begin{mainthm}\label{thm:C}
Let $X$ be a $\sigma$-compact metrizable space, $S$ a finite nonempty
symmetric set of homeomorphisms of $X$, and $m$ a Borel measure on $X$
which is finite on compact sets and whose null sets are preserved by
$S$. Then every $S$-invariant
probability charge on $\Bor(X)$ vanishes on all $m$-null Borel sets if
and only if $(X,S,m)$ has the singular displacement gap.
\end{mainthm}

\Cref{thm:C} combines \cref{thm:null,cor:converse}. The hypothesis
concerns only countably additive measures, while the conclusion concerns
arbitrary finitely additive ones. The two main examples are compact
spaces with probability measures, and Euclidean space with Lebesgue
measure. For actions by isometries
preserving a Borel probability measure of full support (every nonempty open set has
positive measure) on a compact metric space, the singular displacement
gap follows from a spectral gap (\cref{prop:isometric}), and we obtain a
complete characterization.

\begin{mainthm}\label{thm:D}
Let $X$ be a compact metric space, $m$ a Borel probability measure of
full support on $X$, and $S$ a finite nonempty symmetric set of
isometries of $X$ preserving $m$. Then $m$ is the only $S$-invariant
probability charge on $\Bor(X)$ if and only if $(X,S,m)$ has a spectral
gap.
\end{mainthm}

The ``if'' direction is \cref{thm:gap-uniqueness}, and \cref{thm:sphere}
is its special case for spheres. The ``only if'' direction uses a
theorem of Schmidt~\cite[Proposition~2.3]{Schmidt} for nonatomic $m$
(\cref{thm:schmidt}), as explained in the proof of
\cref{thm:isometric-criterion}.

\subsection{Ideas of the proof}
The basic tool is an elementary conservation law
(\cref{lem:certificate}). Regard the action as a graph joining each
point $x$ to the points $sx$, $s\in S$. A \emph{field}
$F=(F_s)_{s\in S}$ of bounded Borel functions prescribes a flow $F_s(x)$
along each edge $x\to sx$, and its \emph{divergence}
\[
 \Div F(x)=\sum_{s\in S}\bigl(F_s(x)-F_s(s^{-1}x)\bigr)
\]
is the net outflow at $x$. An invariant charge $\mu$ integrates every
divergence to zero. Hence if a Borel set $Y$ admits a
\emph{certificate}, a bounded Borel field with $\Div F\ge\ind_Y$ \emph{at
every point}, then every invariant charge vanishes on $Y$.

Bounded Borel flows are also the key tool in the proof of Borel circle
squaring by Marks and Unger~\cite{MarksUnger}.

For the sufficiency in \cref{thm:C} we construct certificates for
invariant null Borel sets $Y$, in two steps. First, every probability
measure on a \emph{compact} null set $K$ is singular, and a Hahn--Banach
argument in the space $C(K)$ of continuous functions on $K$ turns the
singular displacement gap into a nonnegative bounded continuous field
with divergence greater than one on $K$ (\cref{sec:compact}). Second,
and this is the technical core of the paper, these local fields are
combined into a certificate for $Y$: every point selects one of
countably many fixed continuous fields by comparing orbit cutoffs with a
random threshold common to all points, and the selected fields are
averaged over the threshold (\cref{sec:null-theorem}).

For actions by isometries of a compact metric space preserving a Borel
probability measure of full support, a spectral gap implies the singular
displacement gap (\cref{prop:isometric}), so invariant
probability charges vanish on null sets and define invariant functionals
on $L^\infty(m)$. A duality argument based on an $L^1$ form of the
spectral gap shows that such functionals are unchanged when a function is
replaced by its mean (\cref{lem:functional}). This proves
\cref{thm:sphere} and one direction of \cref{thm:D}.

For Euclidean space we apply \cref{thm:null} directly to the group
$\Gamma$ generated by the rotations about $p$ and $q$, with Lebesgue
measure as reference measure. Its singular displacement gap comes from
the spherical spectral gaps and the fact that averaging over spheres
about $p$, $q$, $p$ in turn makes every probability measure absolutely
continuous (\cref{prop:euclid-H}). Continuous $\Gamma$-invariant
functions are constant, so the charge $\mu$ is finite on bounded sets and
grows polynomially on balls (\cref{lem:dense}). Integrating certificates
against $\mu$ over orbit neighbourhoods would force exponential growth
on a null set, so $\mu$ vanishes on bounded null sets
(\cref{lem:euclid-null}). The duality argument on balls then shows that
$\mu$ is unchanged by spherical averaging about $p$ and $q$
(\cref{lem:average-invariance}), three averages and Dini's theorem make
$\mu$ countably additive on bounded sets (\cref{lem:dini}), and a final
mollification identifies $\mu$ with Lebesgue measure.

\subsection{Organization}
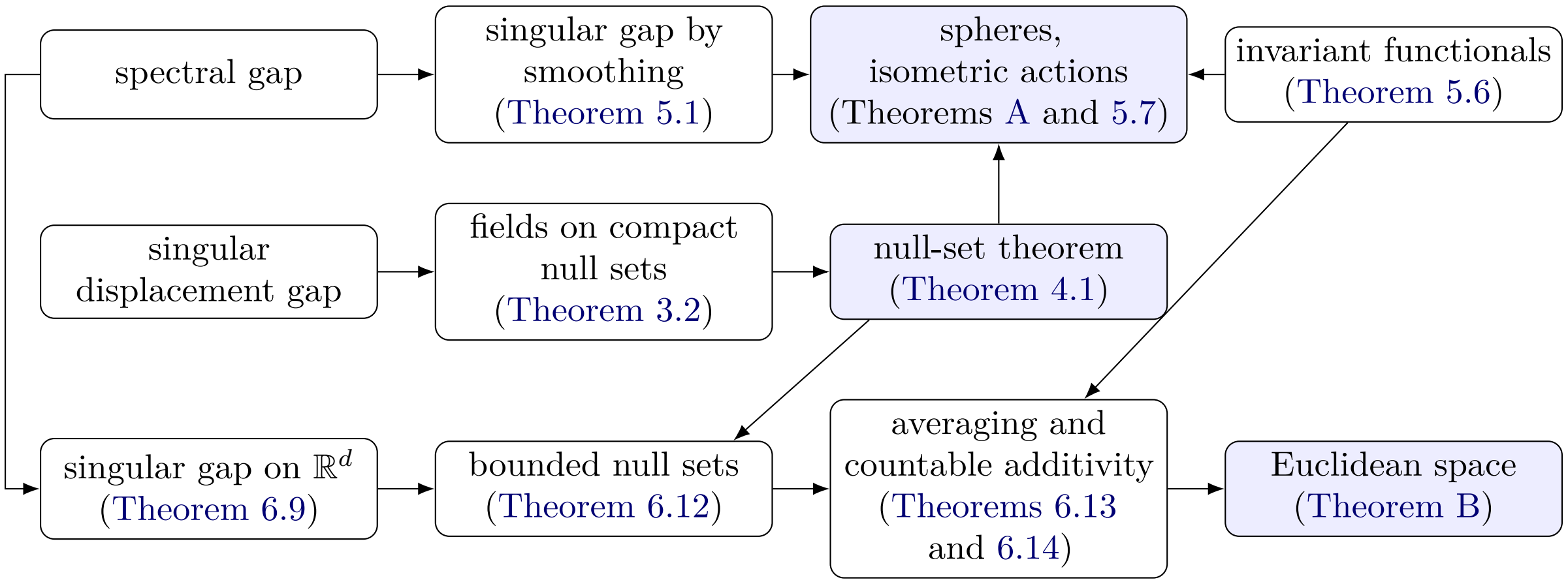
\begin{figure}[tbp]
\centering
\begin{tikzpicture}[>=Latex,font=\small,
 box/.style={draw,rounded corners,align=center,inner sep=3pt,fill=white,minimum height=9mm,text width=3.2cm,font=\small,execute at begin node={\hyphenpenalty=10000}},
 main/.style={box,fill=blue!7}]
 \node[box] (sg) at (0,0) {spectral gap};
 \node[box] (iso) at (4,0) {singular gap by\\smoothing (\cref{prop:isometric})};
 \node[main,text width=3.6cm] (sphere) at (8,0) {spheres, isometric actions\\(Theorems~\ref{thm:sphere} and~\ref{thm:gap-uniqueness})};
 \node[box] (fun) at (12,0) {invariant functionals\\(\cref{lem:functional})};
 \node[box] (H) at (0,-2) {singular\\displacement gap};
 \node[box] (cont) at (4,-2) {fields on compact\\null sets (\cref{lem:continuous-certificate})};
 \node[main] (null) at (8,-2) {null-set theorem\\(\cref{thm:null})};
 \node[box] (rgap) at (0,-4.2) {singular gap on $\R^d$\\(\cref{prop:euclid-H})};
 \node[box] (enull) at (4,-4.2) {bounded null sets\\(\cref{lem:euclid-null})};
 \node[box] (avg) at (8,-4.2) {averaging and\\countable additivity\\(\cref{lem:average-invariance,lem:dini})};
 \node[main] (euclid) at (12,-4.2) {Euclidean space\\(\cref{thm:euclid})};
 \draw[->] (sg)--(iso);
 \draw[->] (iso)--(sphere);
 \draw[->] (fun)--(sphere);
 \draw[->] (H)--(cont);
 \draw[->] (cont)--(null);
 \draw[->] (null)--(sphere);
 \draw[->] ([xshift=4mm]null.south west)--([xshift=-4mm]enull.north east);
 \draw[->] (sg.west)--++(-0.35,0)|-(rgap.west);
 \draw[->] (rgap)--(enull);
 \draw[->] (enull)--(avg);
 \draw[->] (avg)--(euclid);
 \draw[->] (fun)--(avg);
\end{tikzpicture}
\caption{The main proof route. Shaded boxes contain the main results.
The middle row is the general construction. The top row applies it to
isometric actions on compact spaces, and the bottom row to the two
rotation sets acting on $\R^d$ with Lebesgue measure. There the
certificates of the null-set theorem, combined with polynomial growth,
handle bounded null sets, and the Euclidean row ends with mollification.
The spectral gap also enters the invariant-functionals box, through the
$L^1$ form of \cref{lem:L1-gap}. The converse
statements and \cref{thm:D} are omitted from the diagram.}
\label{fig:roadmap}
\end{figure}

\Cref{sec:prelim} introduces charges, fields and the two gap conditions.
\Cref{sec:compact} constructs continuous fields with large divergence
on compact null sets, and \cref{sec:null-theorem} combines them to prove
the null-set theorem, then proves its converse. \Cref{sec:sphere} treats isometric
actions and proves \cref{thm:sphere}, and \cref{sec:euclid} proves
\cref{thm:euclid}. \Cref{sec:complements} proves \cref{thm:D}. \Cref{sec:averages}
collects standard properties of spherical averages. The main
proof route is summarized in \cref{fig:roadmap}.

\phantomsection\label{par:ai}
\subsection*{Use of AI}
The core mathematical ideas and the proofs in this paper were found by
AI models: GPT-6-Astra and GPT-6.1-Sol (OpenAI), and Claude Fable 5.2
(Anthropic). The text was written by these models and Claude Opus 5.5 (Anthropic).
The author's role was organizational and editorial.

\section{Preliminaries}\label{sec:prelim}

This section fixes notation and collects the elementary facts used
throughout. The only non-standard notions are \emph{fields} and their
\emph{divergence} (\cref{def:field}), and the \emph{singular displacement
gap} (\cref{def:H}). A reader familiar with charge integration and total
variation may wish to skim \cref{subsec:conventions} and continue at
\cref{subsec:fields}.

\subsection{Conventions}\label{subsec:conventions}
For a topological space $X$, $\Bor(X)$ is its Borel $\sigma$-algebra. A
\emph{measure} is always countably additive, and a \emph{Borel probability
measure} on $X$ is a countably additive probability measure on $\Bor(X)$.
A Borel measure $\nu$ is \emph{supported on} a Borel set $A$ if
$\nu(X\setminus A)=0$.

Functions are real-valued, and $C(X)$ is the space of continuous
functions $X\to\R$. We write $\ind_A$ for the indicator function of a
set $A$, $\delta_x$ for the Dirac measure at $x$,
$A\triangle B=(A\setminus B)\cup(B\setminus A)$ for the symmetric
difference, and $f^\pm=\max(\pm f,0)$ for the positive and negative parts
of a function. For a measurable space $\Omega$, $\Bb(\Omega)$ is the
space of bounded measurable functions on $\Omega$ with the norm
$\|f\|_\infty=\sup_{\omega\in\Omega}|f(\omega)|$. This is the supremum
over \emph{all} points, not an essential supremum. When a measure $m$ is
fixed, $L^p(m)$ is the usual space of equivalence classes, modulo
equality $m$-almost everywhere, with the norm $\|f\|_{L^p}$, and
$\|f\|_{L^\infty}$ is the essential supremum of $|f|$. We usually do not
distinguish a class from its representatives, except where pointwise
statements matter. For integrable $f$ we write $m(f)=\int f\,dm$. A
finite Borel measure $m$ on a separable metrizable space is
\emph{nonatomic} if $m(\{x\})=0$ for every point $x$.

\subsection{Charges and their integrals}
The objects of study are finitely additive measures on Borel sets. Since we work both with finitely
additive and countably additive measures, to avoid ambiguities we refer to the former as charges.

\begin{definition}[Charges]\label{def:charge}
Let $(\Omega,\mathcal A)$ be a measurable space. A \emph{charge} on
$\Omega$ is a function $\mu\colon\mathcal A\to[0,\infty)$ such that
$\mu(\varnothing)=0$ and $\mu(A\cup B)=\mu(A)+\mu(B)$ whenever
$A,B\in\mathcal A$ are disjoint. It is a \emph{probability charge} if
moreover $\mu(\Omega)=1$. A countably additive charge is a measure. If $s\colon\Omega\to\Omega$ is
an \emph{automorphism} of $(\Omega,\mathcal A)$, that is, a bijection such
that $s$ and $s^{-1}$ are both measurable, then $\mu$ is
\emph{$s$-invariant} if $\mu(sA)=\mu(A)$ for all $A\in\mathcal A$.
\end{definition}

Finite additivity and nonnegativity imply that a charge is monotone and
finitely subadditive. Integration of bounded functions against a charge
is defined by uniform approximation with simple functions. For a simple
function $f=\sum_{i=1}^ka_i\ind_{A_i}$, with $A_1,\dots,A_k\in\mathcal A$
pairwise disjoint, put $\mu(f)=\sum_{i=1}^ka_i\,\mu(A_i)$.
Finite additivity shows that this does not depend on the representation
of $f$, that it is linear and positive on simple functions, and that
$|\mu(f)|\le\mu(\Omega)\|f\|_\infty$. Since simple functions are
uniformly dense in $\Bb(\Omega)$, $\mu$ extends uniquely to a bounded
linear functional on $\Bb(\Omega)$, given by
$\mu(f)=\lim_k\mu(f_k)$ for simple $f_k$ with $\|f_k-f\|_\infty\to0$.
We call $\mu(f)$ the \emph{integral} of $f$ against $\mu$ and use the
same letter for the charge and the functional. The extension is positive
($f\ge0$ implies $\mu(f)\ge0$), satisfies $\mu(\ind_A)=\mu(A)$, and
\begin{equation}\label{eq:integral-bound}
 |\mu(f)|\le\|f\|_\infty\,\mu(\{f\ne0\}) \le \|f\|_\infty\,\mu(\Omega).
\end{equation}
The first inequality follows from positivity and the pointwise bound
$|f|\le\|f\|_\infty\ind_{\{f\ne0\}}$. If $\mu$ is $s$-invariant, then
$\mu(\ind_A\circ s^{-1})=\mu(sA)=\mu(A)$. By linearity and uniform
approximation,
\begin{equation}\label{eq:integral-invariance}
 \mu(f\circ s^{-1})=\mu(f)\qquad(f\in\Bb(\Omega)).
\end{equation}

\subsection{Total variation and singular measures}
For an automorphism $s$ of $\Omega$ and a charge $\nu$, the
\emph{push-forward} $s_*\nu$ is the charge $s_*\nu(A)=\nu(s^{-1}A)$, so
that $s_*\nu(f)=\nu(f\circ s)$. Thus $\nu$ is $s$-invariant exactly
when $s_*\nu=\nu$. If $\lambda$ is a measure on $(\Omega,\mathcal A)$ and
$s_*\lambda=\lambda$, we call $s$ a \emph{measure-preserving automorphism}
of $(\Omega,\lambda)$.

\begin{definition}[Total variation]\label{def:TV}
The \emph{total variation distance} between charges $\mu,\nu$ on
$\Omega$ is
\[
 \|\mu-\nu\|_\TV=\sup\bigl\{|\mu(f)-\nu(f)|:f\in\Bb(\Omega),\ \|f\|_\infty\le1\bigr\}.
\]
\end{definition}

With this normalization two mutually singular probability measures are
at distance two. If $s$ is an automorphism of $\Omega$, then
$f\mapsto f\circ s$ is an isometric bijection of $\Bb(\Omega)$, so
$\|s_*\mu-s_*\nu\|_\TV=\|\mu-\nu\|_\TV$.
For finite Borel measures, the supremum in \cref{def:TV} can be
restricted to continuous test functions in two situations. On locally
compact separable metrizable spaces, such as compact metrizable spaces
and $\R^d$, test functions with compact support suffice
(\cref{lem:TV-continuous}). This is used for the smoothing limits in
\cref{prop:isometric,prop:TV}. For measures supported on a compact set,
bounded continuous test functions suffice on any metrizable space
(\cref{lem:TV-compact}). This is used for the separation argument of
\cref{sec:compact}. For a finite signed Borel measure
$\eta$ we write $|\eta|$ for its total variation measure,
\[
 |\eta|(A)=\sup\Bigl\{\sum_{i=1}^k|\eta(A_i)|:\ A_1,\dots,A_k\ \text{a Borel partition of }A\Bigr\}.
\]
A finite signed measure $\eta$ is \emph{supported on} $A$ if
$|\eta|(X\setminus A)=0$. The following lemma is well known.

\begin{lemma}\label{lem:TV-continuous}
Let $X$ be a locally compact separable metrizable space, and let
$C_c(X)$ be the space of continuous functions on $X$ with compact
support. For every finite signed Borel measure $\eta$ on $X$,
\[
 |\eta|(X)=\sup\Bigl\{\int g\,d\eta:\ g\in C_c(X),\ \|g\|_\infty\le1\Bigr\}.
\]
In particular, if $\mu,\nu$ are finite Borel measures on $X$, then
$\|\mu-\nu\|_\TV=|\mu-\nu|(X)$, and this equals the supremum above for
$\eta=\mu-\nu$.
\end{lemma}

When $X$ is compact, $C_c(X)=C(X)$. The next variant needs no local compactness, only compact support. It is
used in \cref{sec:compact}, where the ambient space need not be locally
compact. For a metrizable space $X$, let $C_b(X)$ be the space of
bounded continuous functions on $X$.

\begin{lemma}\label{lem:TV-compact}
Let $X$ be a metrizable space, and let $\eta$ be a finite signed Borel
measure on $X$ such that $|\eta|$ is supported on a compact set $K$.
Then
\[
 |\eta|(X)=\sup\Bigl\{\int g\,d\eta:\ g\in C_b(X),\ \|g\|_\infty\le1\Bigr\}.
\]
If $\eta=\mu-\nu$ for finite Borel measures $\mu,\nu$, this also equals
$\|\mu-\nu\|_\TV$.
\end{lemma}

\begin{proof}
The inequality $\ge$ is clear. Conversely, $|\eta|(X)=|\eta|(K)$, and by
\cref{lem:TV-continuous} applied to the restriction of $\eta$ to the
compact metrizable space $K$, for every $\eps>0$ there is $g_0\in C(K)$
with $\|g_0\|_\infty\le1$ and $\int_Kg_0\,d\eta>|\eta|(X)-\eps$. By the
Tietze extension theorem, $g_0$ extends to
some $g\in C_b(X)$ with $\|g\|_\infty\le1$, and $\int g\,d\eta=\int_Kg_0\,d\eta$.
The last assertion follows from a Hahn decomposition, which uses no
topology.
\end{proof}

\begin{definition}[Singular measures]\label{def:singular}
Let $m$ be a Borel measure on $X$. A Borel probability
measure $\nu$ is \emph{singular} with respect to $m$, written
$\nu\perp m$, if it is supported on an $m$-null Borel set:
\[
 \nu\perp m\iff\exists N\in\Bor(X)\quad m(N)=0\ \text{and}\ \nu(X\setminus N)=0.
\]
\end{definition}

\subsection{Actions, fields and divergence}\label{subsec:fields}
In the rest of this section and in \cref{sec:compact,sec:null-theorem},
$X$ is a \emph{$\sigma$-compact} metrizable space, that is, a countable
union of compact sets, $S$ is a finite nonempty set of homeomorphisms of
$X$ which is \emph{symmetric}, $S=S^{-1}$, and $\Gamma=\langle S\rangle$ is
the countable group it generates.
We use the word length and its balls, for integers $j\ge0$,
\[
 |g|_S=\min\{k\ge0:g=s_1\cdots s_k,\ s_i\in S\},\qquad
 \Gamma_j=\{g\in\Gamma:|g|_S\le j\},
\]
where the empty word ($k=0$) represents the identity.
\begin{definition}[Invariance]\label{def:invariance}
A Borel measure $m$ on $X$ is \emph{quasi-invariant} if $\Gamma$
preserves its null sets: $m(A)=0$ if and only if $m(gA)=0$, for all
$A\in\Bor(X)$ and $g\in\Gamma$. It is \emph{invariant} if $m(sA)=m(A)$
for all $A\in\Bor(X)$ and $s\in S$, and it has \emph{full support} if
$m(U)>0$ for every nonempty open $U\subseteq X$. The \emph{saturation}
of a set $A\subseteq X$ is $\Gamma A=\bigcup_{g\in\Gamma}gA$, and $A$ is
\emph{invariant} if $\Gamma A=A$. A probability charge $\mu$ on
$\Bor(X)$ is \emph{$S$-invariant} if it is $s$-invariant for every
$s\in S$.
\end{definition}

\begin{figure}[ht]
\centering
\begin{tikzpicture}[>=Latex,font=\small,
 pt/.style={circle,fill,inner sep=1.4pt},
 every edge/.style={draw},shorten >=2pt,shorten <=2pt]
 \node[pt,label={[label distance=4pt,fill=white,inner sep=1pt]-135:$x$}] (x) at (0,0) {};
 \node[pt,label=right:$ax$] (ax) at (3.2,0) {};
 \node[pt,label=left:$a^{-1}x$] (aix) at (-3.2,0) {};
 \node[pt,label=above:$bx$] (bx) at (0,2.2) {};
 \node[pt,label=below:$b^{-1}x$] (bix) at (0,-2.2) {};
 \draw[->,thick,blue!60!black] (x) to[bend left=14] node[above,pos=.55] {$F_a(x)$} (ax);
 \draw[->,thick,red!60!black] (ax) to[bend left=14] node[below,pos=.45] {$F_{a^{-1}}(ax)$} (x);
 \draw[->,thick,blue!60!black] (x) to[bend left=14] node[below,pos=.55] {$F_{a^{-1}}(x)$} (aix);
 \draw[->,thick,red!60!black] (aix) to[bend left=14] node[above,pos=.45] {$F_a(a^{-1}x)$} (x);
 \draw[->,thick,blue!60!black] (x) to[bend left=14] node[left,pos=.55] {$F_b(x)$} (bx);
 \draw[->,thick,red!60!black] (bx) to[bend left=14] node[right,pos=.45] {$F_{b^{-1}}(bx)$} (x);
 \draw[->,thick,blue!60!black] (x) to[bend left=14] node[right,pos=.55] {$F_{b^{-1}}(x)$} (bix);
 \draw[->,thick,red!60!black] (bix) to[bend left=14] node[left,pos=.45] {$F_b(b^{-1}x)$} (x);
\end{tikzpicture}
\caption{The Schreier graph near $x$ for $S=\{a^{\pm1},b^{\pm1}\}$. The
picture is schematic: the five points need not be distinct. The field
sends $F_s(x)$ from $x$ to $sx$ (blue arrows, outgoing terms) and
receives $F_s(s^{-1}x)$ from $s^{-1}x$ (red arrows, incoming terms). The
values themselves may have either sign. The divergence $\Div F(x)$ is
the total outflow minus the total inflow. A certificate for $Y$ makes
every point of $Y$ a net source of at least one unit, and no point a net
sink.}
\label{fig:divergence}
\end{figure}
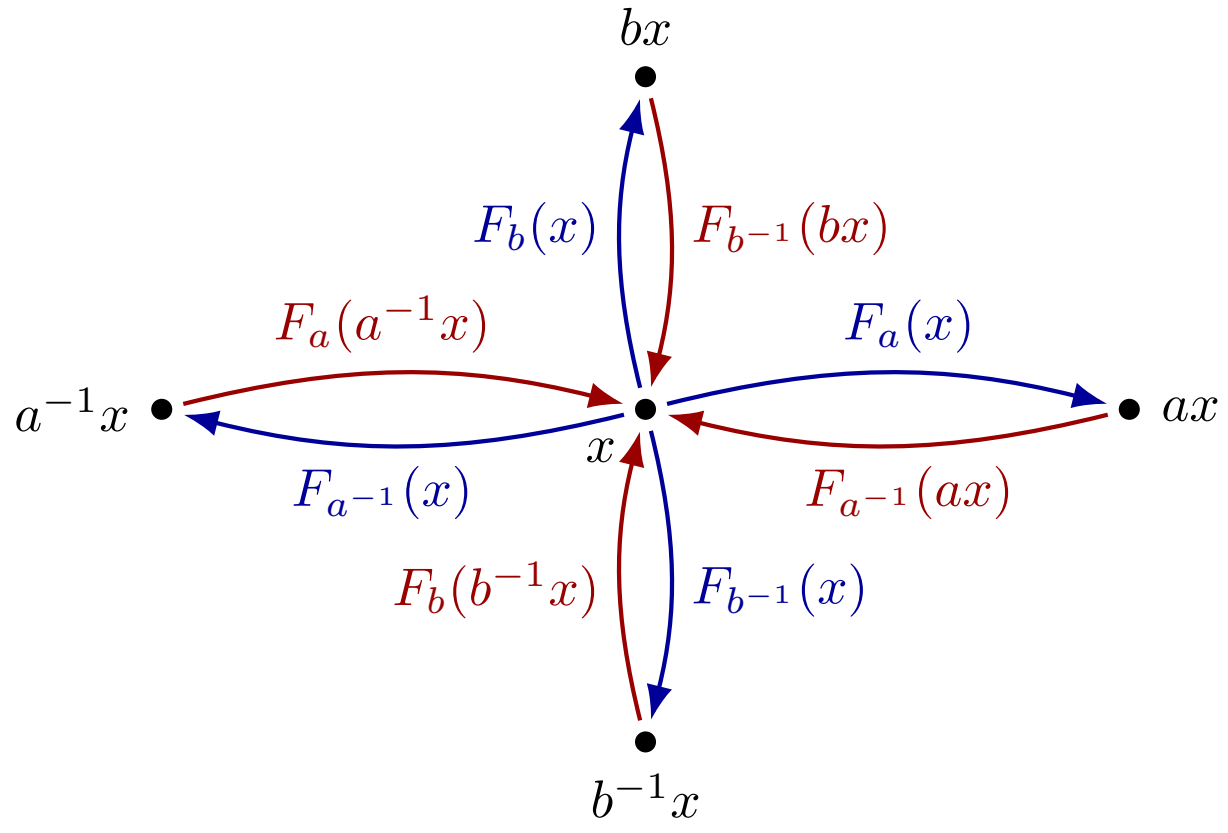

The following notion is the basic tool of the paper. Think of the
Schreier graph of the action, which joins each point $x$ to the points
$sx$, $s\in S$. A field prescribes a real ``flow'' $F_s(x)$ along each
oriented edge $x\to sx$, and its divergence at $x$ is the net amount
leaving $x$ (see \cref{fig:divergence}). The edges are labelled by
generators. Loops ($sx=x$) and multiple edges are allowed, since the
action need not be free. Unlike a flow in the sense of graph theory, a
field is not required to satisfy any antisymmetry condition such as
$F_{s^{-1}}(sx)=-F_s(x)$.

\begin{definition}[Fields and divergence]\label{def:field}
A \emph{field} on $X$ is a family of bounded Borel functions indexed by
$S$, and its \emph{sum norm} is the sum of their norms:
\[
 F=(F_s)_{s\in S}\in\Bb(X)^S,\qquad \|F\|_1=\sum_{s\in S}\|F_s\|_\infty .
\]
Its \emph{divergence} is the function
\[
 \Div F(x)=\sum_{s\in S}\Bigl(F_s(x)-F_s(s^{-1}x)\Bigr)
 \qquad(x\in X).
\]
The field is \emph{nonnegative} if every $F_s\ge0$, and it is
\emph{supported on} a set $Y\subseteq X$ if every component vanishes
outside $Y$.
A \emph{certificate} for a Borel set $Y\subseteq X$ is a field whose
divergence dominates the indicator of $Y$ at every point:
\[
 \Div F(x)\ge\ind_Y(x)\qquad(x\in X).
\]
\end{definition}

Invariant charges ``conserve'' divergence, and this is how certificates
are used.

\begin{lemma}[Certificates annihilate sets]\label{lem:certificate}
Let $\mu$ be an $S$-invariant probability charge on $\Bor(X)$. Then
$\mu(\Div F)=0$ for every field $F$. Consequently, if a Borel set $Y$
admits a certificate, then $\mu(Y)=0$.
\end{lemma}

\begin{proof}
For each $s$, $\mu(F_s\circ s^{-1})=\mu(F_s)$ by
\eqref{eq:integral-invariance}. Summing over $s$ gives $\mu(\Div F)=0$.
If $\Div F\ge\ind_Y$, positivity of $\mu$ gives
$\mu(Y)\le\mu(\Div F)=0$.
\end{proof}

The proof uses nothing about $\mu$ beyond finite additivity and
invariance. In particular, it applies to charges which do not vanish on
null sets. The difficulty, of course, is to construct certificates. Since
the integral against a charge can depend on the values of a function on
a set of reference measure zero, both the bound on $F$ and the
inequality $\Div F\ge\ind_Y$ must hold at \emph{every} point, not merely
almost everywhere.

\subsection{The two gap conditions}
Throughout, the reference measure $m$ is a Borel measure finite on
compact sets. It need not be a probability measure. For instance, $m$
may be Lebesgue measure on $\R^d$.

\begin{definition}[Singular displacement gap]\label{def:H}
Let $m$ be a Borel measure on $X$ and let $\delta>0$. The
system $(X,S,m)$ satisfies the \emph{singular displacement gap} with
constant $\delta$ if
\begin{equation*}\tag{H}\label{eq:H}
 \max_{s\in S}\|s_*\nu-\nu\|_\TV\ge\delta
 \qquad\text{for every Borel probability measure }\nu\perp m.
\end{equation*}
\end{definition}

For instance, \eqref{eq:H} fails if some $\nu\perp m$ is $S$-invariant.
This happens when all elements of $S$ fix a point $z$ with
$m(\{z\})=0$, since then $\nu=\delta_z$ is invariant and singular.
Condition \eqref{eq:H} is necessary for invariant charges to vanish on
null sets (\cref{cor:converse}). For isometries of a compact metric
space preserving a probability measure of full support, it follows from
the classical spectral gap condition (\cref{prop:isometric}).

\begin{definition}[Spectral gap]\label{def:spectral-gap}
Let $m$ be an $S$-invariant Borel probability measure on $X$ and let
$\kappa>0$. The system $(X,S,m)$ has a \emph{spectral gap} with constant
$\kappa$ if
\begin{equation}\label{eq:spectral-gap}
 \max_{s\in S}\|u\circ s^{-1}-u\|_{L^2(m)}\ge\kappa\,\|u-m(u)\|_{L^2(m)}
 \qquad\bigl(u\in L^2(m)\bigr).
\end{equation}
\end{definition}

\begin{definition}[Ergodicity]\label{def:ergodic}
Let $m$ be an $S$-invariant Borel probability measure on $X$. The action
is \emph{ergodic} if every Borel set which is invariant up to null sets
is null or conull:
\[
 m(sA\,\triangle\,A)=0\ \text{for all }s\in S\implies m(A)\in\{0,1\}\qquad(A\in\Bor(X)).
\]
\end{definition}

Taking $u=\ind_A$ in \eqref{eq:spectral-gap} shows that a spectral gap
implies ergodicity.

\section{Fields with large divergence on compact null sets}\label{sec:compact}

In this section $X$ is a $\sigma$-compact metrizable space, $S$ a finite
nonempty symmetric set of homeomorphisms, and $m$ a Borel measure on $X$
satisfying \eqref{eq:H} with constant $\delta>0$. Quasi-invariance is
not needed here.

Every probability measure supported on a compact $m$-null set $K$ is singular with respect to $m$,
so by \eqref{eq:H} it is moved by some generator. \Cref{lem:continuous-certificate} shows the
existence of a field of controlled size whose divergence exceeds one on $K$.
In \cref{sec:null-theorem} the resulting fields are combined into
certificates for arbitrary invariant null sets.

\begin{definition}[Continuous fields]\label{def:continuous-field}
A field $F$ on $X$ is \emph{continuous} if all its components are bounded
continuous functions:
\[
 F\in C_b(X)^S .
\]
\end{definition}

\begin{lemma}[Continuous fields on a compact null set]\label{lem:continuous-certificate}
Let $K\subseteq X$ be a compact set with $m(K)=0$. There is a
nonnegative continuous field $c$ such that
\[
 \|c\|_1<4/\delta\qquad\text{and}\qquad\Div c(x)>1\quad(x\in K),
\]
where $\|c\|_1=\sum_{s\in S}\|c_s\|_\infty$.
\end{lemma}

\begin{proof}
If $K=\varnothing$ take $c=0$. Otherwise put $C=2/\delta$ and consider
the following two convex subsets of the Banach space $C(K)$:
\[
 V=\bigl\{(\Div F)|_K:\ F\in C_b(X)^S,\ \|F\|_1<C\bigr\},\qquad
 U=\bigl\{u\in C(K):\ u(x)>1\text{ for all }x\in K\bigr\}.
\]
Here $(\Div F)|_K$ is the restriction to $K$ of the continuous function
$\Div F$ on $X$. Its values on $K$ involve the values of $F$ on
$K\cup\bigcup_{s\in S}s^{-1}K$. We show that $U\cap V\ne\varnothing$.

Suppose instead that $U\cap V=\varnothing$. Since $U$ is open and
nonempty, the separation theorem for convex sets
 gives a continuous linear functional
$\ell$ on $C(K)$ and $\gamma\in\R$ with
\[
 \ell(v)\le\gamma<\ell(u)\qquad(v\in V,\ u\in U).
\]
Taking $v=0\in V$ shows $\gamma\ge0$. If $h\in C(K)$ is nonnegative, then
$2+ah\in U$ for all $a\ge0$, so $0\le\gamma<\ell(2)+a\,\ell(h)$ for all
$a\ge0$. This forces $\ell(h)\ge0$. Thus $\ell$ is a positive functional.
It is nonzero, so $\ell(1)>0$, and after dividing $\ell$ and $\gamma$ by
$\ell(1)$ we may assume $\ell(1)=1$. By the Riesz representation theorem
for positive functionals on $C(K)$ there is
a Borel probability measure $\nu_K$ on $K$ representing $\ell$, and we
extend it to $X$:
\[
 \ell(u)=\int_Ku\,d\nu_K\quad(u\in C(K)),\qquad
 \nu(A)=\nu_K(A\cap K)\quad(A\in\Bor(X)).
\]
Then $\nu$ is supported on $K$ and $m(K)=0$, so $\nu\perp m$. Testing
against the constants $1+\eps\in U$ gives $\gamma\le1$.

Now fix $s_0\in S$ and $g\in C_b(X)$ with $\|g\|_\infty\le1$, let
$0<a<C = 2/\delta$, and consider the field $F$ with $F_{s_0}=ag$ and $F_s=0$ for
$s\ne s_0$. Then $(\Div F)|_K\in V$, so
\[
 a\Bigl(\int g\,d\nu-\int g\,d(s_0^{-1})_*\nu\Bigr)
 =\int_K\Div F\,d\nu=\ell\bigl((\Div F)|_K\bigr)\le\gamma\le1.
\]
The signed measure $\nu-(s_0^{-1})_*\nu$ is supported on the compact set
$K\cup s_0^{-1}K$. Taking the supremum over $g$ and $a$, and using
\cref{lem:TV-compact},
\[
 C\,\|(s_0)_*\nu-\nu\|_\TV=C\,\|\nu-(s_0^{-1})_*\nu\|_\TV\le1 .
\]
The equality holds because $(s_0)_*$ preserves total variation. Since
$s_0\in S$ was arbitrary, $C\max_s\|s_*\nu-\nu\|_\TV\le1<C\delta$,
contradicting \eqref{eq:H}.

Therefore some $F\in C_b(X)^S$ with $\|F\|_1<C$ satisfies $\Div F>1$ on
$K$. Shift each component by its norm:
\[
 c_s=F_s+\|F_s\|_\infty\qquad(s\in S).
\]
Adding a constant to a component does not change the divergence, since
the constant cancels in $c_s-c_s\circ s^{-1}$. Thus $\Div c=\Div F>1$ on
$K$, $c_s\ge0$, and $\|c\|_1\le2\|F\|_1<2C=4/\delta$.
\end{proof}

Only the divergence on $K$ is controlled. Off $K$ it may be negative, so
the construction does not in general give a certificate for $K$ in the
sense of \cref{def:field}. The construction in \cref{sec:null-theorem} is
designed so that the inequality $\Div c>1$ on a single compact set is
all that is ever used. Since $\Div c$ is continuous, the set
\[
 O=\{x\in X:\Div c(x)>1\}
\]
is an open neighbourhood of $K$, and the same field works for every
compact set contained in $O$.

\section{The null-set theorem}\label{sec:null-theorem}

The main obstacle in the Borel Ruziewicz problem is that an invariant
probability charge on Borel sets might give positive mass to a Borel set
of reference measure zero. The following theorem removes this obstacle
under the singular displacement gap.

\Needspace{10\baselineskip}
\begin{theorem}[Null-set theorem]\label{thm:null}
Let $X$ be a $\sigma$-compact metrizable space, $S$ a finite nonempty
symmetric set of homeomorphisms of $X$, and $m$ a quasi-invariant Borel
measure on $X$ which is finite on compact sets and satisfies \eqref{eq:H}
with constant $\delta>0$. Then every invariant $m$-null Borel set
$Y\subseteq X$ admits a nonnegative certificate supported on $Y$ with
components bounded by $8/\delta$: there is a field $H$ on $X$ such that
\[
 0\le H_s\le 8/\delta,\qquad H_s=0\ \text{off}\ Y\quad(s\in S),
 \qquad \Div H\ge\ind_Y .
\]
Consequently $\mu(N)=0$ for every $S$-invariant probability charge $\mu$
on $\Bor(X)$ and every $m$-null Borel set $N\subseteq X$.
\end{theorem}

The last sentence follows from the main part of the theorem. If $N$ is an $m$-null Borel
set, its saturation $Y=\Gamma N$ is an invariant $m$-null Borel set, by
quasi-invariance. By \cref{lem:certificate}, $\mu(N)\le\mu(Y)=0$. The
inequality $\Div H\ge\ind_Y$ holds at every point of $X$.

The two main cases are a compact metrizable space with a probability
measure (used for spheres in \cref{sec:sphere}) and Euclidean space with
Lebesgue measure (used in \cref{sec:euclid}). The measure $m$ enters only through
its restrictions to compact sets and their finitely many translates.

In \cref{subsec:covers,subsec:cutoffs} we
construct compact covers of $Y$ and orbit cutoffs. In
\cref{subsec:selection} a threshold common to all points determines a
selection among fixed continuous fields. Averaging over the threshold
gives the certificate in \cref{subsec:averaged}. We prove the converse in
\cref{subsec:converse}.

Throughout \cref{subsec:covers,subsec:cutoffs,subsec:selection,subsec:averaged}
we fix $X$, $S$, $m$, $\delta$ and the invariant $m$-null Borel set $Y$
as in \cref{thm:null}, and a compatible metric $\varrho$ on $X$. We may
assume $X\ne\varnothing$. Put $M=4/\delta$.

\subsection{Compact covers}\label{subsec:covers}
A null set need not lie in any compact null set. Instead we cover $Y$ by compact sets
of prescribed small measures, in such a way that every point of $Y$ is covered infinitely
often. Then every point of $Y$ lies in a compact intersection of covering sets whose measures tend
to zero, and such an intersection is a compact null set. The $r$th covering set will later be
translated by many group elements, so its smallness is measured by a different measure $m_r$ for
each $r$.

\begin{lemma}[Compact limsup covers]\label{lem:covers}
Let $(m_r)_{r\ge1}$ be Borel measures on $X$ which are finite on
compact sets and vanish on $Y$, and let $(\eps_r)_{r\ge1}$ be positive
reals. There are compact sets $C_r\subseteq X$ such that
\[
 m_r(C_r)\le\eps_r\quad(r\ge1)\qquad\text{and}\qquad
 Y\subseteq\limsup_{r\to\infty}C_r=\bigcap_{k\ge1}\bigcup_{r\ge k}C_r .
\]
In words, every point of $Y$ lies in $C_r$ for infinitely many $r$.
\end{lemma}

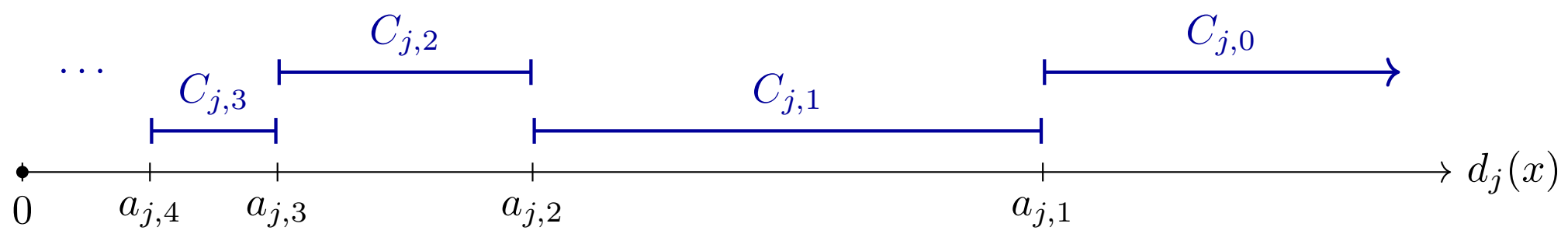
\begin{figure}[ht]
\centering
\begin{tikzpicture}[x=1.45cm,font=\small]
 \draw[->] (0,0)--(8.4,0) node[right] {$d_j(x)$};
 \foreach \x/\l in {0/0,0.75/a_{j,4},1.5/a_{j,3},3/a_{j,2},6/a_{j,1}}
   {\draw (\x,.08)--(\x,-.08) node[below] {$\l$};}
 \draw[|-|,thick,blue!60!black] (3,.35)--(6,.35) node[midway,above] {$C_{j,1}$};
 \draw[|-|,thick,blue!60!black] (1.5,.85)--(3,.85) node[midway,above] {$C_{j,2}$};
 \draw[|-|,thick,blue!60!black] (.75,.35)--(1.5,.35) node[midway,above] {$C_{j,3}$};
 \draw[|->,thick,blue!60!black] (6,.85)--(8.1,.85) node[midway,above] {$C_{j,0}$};
 \node[blue!60!black] at (.37,.85) {$\cdots$};
 \fill[black] (0,0) circle (1.5pt);
\end{tikzpicture}
\caption{The closed bands $C_{j,i}$, drawn through the corresponding
intervals of the distance $d_j$ to $X_j\setminus U_j$. Consecutive
distance intervals share an endpoint. Since $a_{j,i}\downarrow0$, the
bands cover $U_j=\{d_j>0\}$. The endpoints $a_{j,i}$ are chosen adaptively, so that
the $i$th band is small for its own measure $m_{j,i}$.}
\label{fig:bands}
\end{figure}

\begin{proof}
Write $X=\bigcup_jX_j$ with compact sets $X_1\subseteq X_2\subseteq\cdots$,
and fix a bijection $r\leftrightarrow(j,i)$ between the positive integers
and the pairs with $j\ge1$, $i\ge0$. We write $m_{j,i}$, $\eps_{j,i}$
and $C_{j,i}$ for the objects with index $r\leftrightarrow(j,i)$.

Fix $j$. The restriction of $m_{j,0}$ to the compact metrizable space
$X_j$ is a finite Borel measure, hence outer regular:
\[
 m_{j,0}(A)=\inf\{m_{j,0}(U):U\supseteq A,\ U\text{ relatively open in }X_j\}\qquad(A\in\Bor(X_j)).
\]
Since $m_{j,0}(Y)=0$, there is a relatively open set $U_j\subseteq X_j$
with $Y\cap X_j\subseteq U_j$ and $m_{j,0}(U_j)<\eps_{j,0}$. Measure
distances to its complement in $X_j$:
\[
 d_j(x)=\varrho(x,X_j\setminus U_j)=\inf_{y\in X_j\setminus U_j}\varrho(x,y)\qquad(x\in X_j),
\]
with $d_j=\infty$ if $U_j=X_j$. Then $U_j=\{x\in X_j:d_j(x)>0\}$. As
$a\downarrow0$, the sets $\{x\in X_j:0<d_j(x)\le a\}$ decrease to the
empty set, and each $m_{j,i}$ is finite on $X_j$. So we can choose
reals $a_{j,1}>a_{j,2}>\cdots>0$ tending to zero such that
\[
 m_{j,i}\bigl\{x\in X_j:0<d_j(x)\le a_{j,i}\bigr\}\le\eps_{j,i}\qquad(i\ge1).
\]
Cut $U_j$ into the closed bands
\[
 C_{j,0}=\{x\in X_j:d_j(x)\ge a_{j,1}\},\qquad
 C_{j,i}=\{x\in X_j:a_{j,i+1}\le d_j(x)\le a_{j,i}\}\quad(i\ge1),
\]
see \cref{fig:bands}. They are closed subsets of $X_j$, hence compact.
The outer band lies in $U_j$, so $m_{j,0}(C_{j,0})<\eps_{j,0}$. Each
inner band has $m_{j,i}(C_{j,i})\le\eps_{j,i}$ by the choice of
$a_{j,i}$. The bands cover $U_j$, because $a_{j,i}\to0$.

Every $x\in Y$ lies in $X_j$, hence in $U_j$, for all large $j$. So it
lies in at least one band $C_{j,i}$ for each large $j$, and therefore in
$C_r$ for infinitely many $r$.
\end{proof}

\FloatBarrier

\subsection{Orbit cutoffs}\label{subsec:cutoffs}
We will select fields by comparing orbit cutoffs with one threshold
$\omega\in[0,1)$, used at every point. If the $r$th cutoff changes by at
most $1/L_r$ along an edge, the corresponding comparisons at two
neighbouring points disagree on a set of thresholds of measure at most
$1/L_r$. The selections will agree whenever all comparisons agree. Since
each field component lies in $[0,M)$, averaging over the threshold loses
at most $M|S|\sum_{r\ge1}1/L_r$ in divergence. To retain divergence at
least $1/2$, we therefore fix integers $L_r\ge1$ with
\begin{equation}\label{eq:budget}
 M\,|S|\sum_{r\ge1}\frac1{L_r}\le\frac12 ,
\end{equation}
for instance $L_r=2^r\lceil 2M|S|\rceil$. The $r$th covering set will be
thickened by $L_r$ steps along orbits, so we measure it by
\[
 m_r(A)=\sum_{g\in\Gamma_{L_r}}m(gA)\qquad(A\in\Bor(X)).
\]
Recall that $\Gamma_k=\{g\in\Gamma:|g|_S\le k\}$ is the ball of radius
$k$ in the word metric of $\Gamma=\langle S\rangle$.
Each $m_r$ is a Borel measure which is finite on compact sets, since
$gK$ is compact for compact $K$, and $m_r(Y)=0$ by quasi-invariance.
Let $(C_r)$ be compact sets as in \cref{lem:covers}, with these measures
and with $\eps_r=2^{-r}$.

\begin{definition}[Thickened covers and orbit cutoffs]\label{def:cutoff}
The $r$th \emph{thickened cover} $E_r$ and the $r$th \emph{orbit cutoff}
$\phi_r$ are the union of the translates of $C_r$ by words of length at
most $L_r$, and the average of the indicators of the intermediate
thickenings:
\begin{equation}\label{eq:cutoffs}
 E_r=\Gamma_{L_r}C_r=\bigcup_{g\in\Gamma_{L_r}}gC_r,\qquad
 \phi_r(x)=\frac1{L_r}\sum_{k=0}^{L_r-1}\ind_{\Gamma_kC_r}(x).
\end{equation}
\end{definition}

Each $E_r$ is compact, being a finite union of compact sets, and
\begin{equation}\label{eq:Er-small}
 m(E_r)\le\sum_{g\in\Gamma_{L_r}}m(gC_r)=m_r(C_r)\le2^{-r}.
\end{equation}
The orbit cutoff $\phi_r$ is a Borel function with values in $[0,1]$.

\begin{lemma}\label{lem:cutoff}
For all $x\in X$, $s\in S$ and $r\ge1$:
\begin{enumerate}[label=\textup{(\alph*)}]
\item\label{it:cut-pos} If $\phi_r(x)>0$, then $x\in E_r$.
\item\label{it:cut-full} If $x\in C_r$, then $\phi_r(x)=1$.
\item\label{it:cut-lip} $|\phi_r(sx)-\phi_r(x)|\le1/L_r$.
\end{enumerate}
\end{lemma}

\begin{proof}
\ref{it:cut-pos} If $\phi_r(x)>0$, then $x\in\Gamma_kC_r\subseteq E_r$
for some $k<L_r$. \ref{it:cut-full} The identity lies in every
$\Gamma_k$. \ref{it:cut-lip} If $sx\in\Gamma_kC_r$, then
$x\in s^{-1}\Gamma_kC_r\subseteq\Gamma_{k+1}C_r$. Hence, as the sets
$\Gamma_kC_r$ increase with $k$,
\[
 L_r\phi_r(sx)\le\sum_{k=1}^{L_r}\ind_{\Gamma_kC_r}(x)\le L_r\phi_r(x)+1 .
\]
Applying this to the point $sx$ and the generator $s^{-1}\in S$ gives the
reverse inequality.
\end{proof}

\subsection{Fixed fields and selection}\label{subsec:selection}
The compact sets used at a point will be finite intersections of the
sets $E_r$. Such an intersection need not be null, so
\cref{lem:continuous-certificate} does not apply to it directly. We only
use those finite intersections on which a suitable field exists.

\begin{definition}[Admissible index sets]\label{def:admissible}
For a nonempty finite set $F$ of positive integers put
\[
 E_F=\bigcap_{r\in F}E_r ,
\]
a compact set. The set $F$ is \emph{admissible} if some nonnegative
continuous field works on $E_F$ with the bound $M=4/\delta$ from
\cref{lem:continuous-certificate}:
\[
 F\text{ is admissible}\iff\exists\,c\in C_b(X)^S\quad
 \bigl(c_s\ge0\ (s\in S),\ \ \|c\|_1<M,\ \ \Div c>1\text{ on }E_F\bigr).
\]
\end{definition}

There are only countably many finite sets $F$. For each admissible $F$
we choose, once and for all, one such field $c^F$.

\begin{definition}[Bits and selected field]\label{def:selection}
For a \emph{threshold} $\omega\in[0,1)$, the $r$th \emph{bit} at $x$, \(b_{r}(x,\omega)\),
records whether the cutoff exceeds the threshold, and the $n$th
\emph{active prefix}, \(F_{n}(x,\omega)\), collects the active indices up to $n$:
\[
 b_r(x,\omega)=\begin{cases}1,&\omega<\phi_r(x),\\0,&\text{otherwise,}\end{cases}\qquad
 F_n(x,\omega)=\{r\le n:\ b_r(x,\omega)=1\}.
\]
The \emph{selection index} is the first $n$ for which the active prefix
is nonempty and admissible, and the \emph{selected field} is the
corresponding fixed field:
\[
 n(x,\omega)=\min\bigl\{n\ge1:\ F_n(x,\omega)\ne\varnothing\ \text{and}\ F_n(x,\omega)\text{ is admissible}\bigr\},
\]
with the convention $\min\varnothing=\infty$, and
\[
 c(x,\omega)=\begin{cases}c^{F_{n(x,\omega)}(x,\omega)},&n(x,\omega)<\infty,\\
 0,&n(x,\omega)=\infty.\end{cases}
\]
\end{definition}

Note that $c(x,\omega)$ is a field on $X$, and $c_s(x,\omega)(z)$ denotes
the value of its component $c_s(x,\omega)$ at a point $z\in X$. The
next lemma gives the measurability needed for threshold averaging and
the divergence calculation.

\begin{lemma}\label{lem:field-measurable}
For each $s\in S$, the function $(x,\omega,z)\mapsto c_s(x,\omega)(z)$ is
Borel on $X\times[0,1)\times X$, with values in $[0,M)$. Consequently
the functions
\[
 (x,\omega)\mapsto c_s(x,\omega)(x),\qquad(x,\omega)\mapsto c_s(x,\omega)(s^{-1}x),\qquad
 (x,\omega)\mapsto c_s(s^{-1}x,\omega)(s^{-1}x)
\]
are Borel.
\end{lemma}

\begin{proof}
Each bit is a Borel function on $X\times[0,1)$, since $\phi_r$ is Borel.
For $n\ge1$ and a nonempty admissible $F\subseteq\{1,\dots,n\}$, the event
that $n(x,\omega)=n$ and $F_n(x,\omega)=F$ depends only on the first $n$
bits. It is a finite Boolean combination of the Borel sets
$\{b_r=1\}$, $r\le n$, so it is Borel. On it, $c_s(x,\omega)(z)=c^F_s(z)$,
a continuous function of $z$. On the remaining Borel set, where
$n(x,\omega)=\infty$, $c_s(x,\omega)(z)=0$. These countably many
Borel pieces give the claim. The values lie in $[0,M)$ because every
fixed field is nonnegative with $\|c^F\|_1<M$. The three displayed
functions are compositions of this function with the Borel maps
$(x,\omega)\mapsto(x,\omega,x)$, $(x,\omega)\mapsto(x,\omega,s^{-1}x)$ and
$(x,\omega)\mapsto(s^{-1}x,\omega,s^{-1}x)$.
\end{proof}

The next lemma contains everything the divergence estimate needs.
Its two ingredients are that the selection succeeds at points of $Y$,
and that neighbouring selections disagree only for a small set of
thresholds. Let $\lambda$ denote Lebesgue measure on the threshold space:
\[
 \lambda(A)=\int_0^1\ind_A(\omega)\,d\omega\qquad\bigl(A\in\Bor([0,1))\bigr).
\]

\begin{lemma}\label{lem:selected}
Let $x\in X$.
\begin{enumerate}[label=\textup{(\alph*)}]
\item\label{it:sel-null} If $x\in Y$, then for every $\omega$ the
selected field $c=c(x,\omega)$ is one of the fixed fields, and
$\Div c(x)>1$.
\item\label{it:sel-couple} For every $s\in S$,
\[
 \lambda\bigl\{\omega\in[0,1):\ b_r(x,\omega)\ne b_r(s^{-1}x,\omega)\text{ for some }r\bigr\}
 \le\sum_{r\ge1}\frac1{L_r}.
\]
If all bits of $x$ and $s^{-1}x$ agree at $\omega$, then
$c(x,\omega)=c(s^{-1}x,\omega)$.
\end{enumerate}
\end{lemma}

\begin{proof}
\ref{it:sel-null} By \cref{lem:covers} there are infinitely many $r$ with
$x\in C_r$. For these, $\phi_r(x)=1>\omega$ by
\cref{lem:cutoff}\ref{it:cut-full}, so $b_r(x,\omega)=1$. Let $n_0$ be the
least active index. The active prefixes $F_n=F_n(x,\omega)$ are nonempty
for $n\ge n_0$, and the compact sets $E_{F_n}$, $n\ge n_0$, decrease.
Their intersection
\[
 K=\bigcap\bigl\{E_r:\ b_r(x,\omega)=1\bigr\}
\]
contains $x$ by \cref{lem:cutoff}\ref{it:cut-pos}, and it is contained in
$E_r$ for infinitely many $r$. By \eqref{eq:Er-small}, $K$ is a compact
$m$-null set. \Cref{lem:continuous-certificate} gives a nonnegative
continuous field $c'$ with $\|c'\|_1<M$ and $\Div c'>1$ on $K$. The set
$O=\{\Div c'>1\}$ is open and contains $K$. The compact sets
$E_{F_n}\setminus O$, $n\ge n_0$, decrease and have empty intersection,
so one of them is empty, that is, $E_{F_n}\subseteq O$ for some $n\ge n_0$. Then $F_n$ is
admissible, so $n(x,\omega)<\infty$. The selected field is $c^F$ with
$F=F_{n(x,\omega)}(x,\omega)$. Every $r\in F$ is active at $x$, so $x\in
E_F$ by \cref{lem:cutoff}\ref{it:cut-pos}, and $\Div c^F(x)>1$ by
admissibility of $F$.

\ref{it:sel-couple} Fix $r$. For each point $y$, the set of thresholds
activating the $r$th bit is the interval $[0,\phi_r(y))$. The bits at
$x$ and $s^{-1}x$ differ exactly on the symmetric difference of two such
intervals. Its length is at most $|\phi_r(x)-\phi_r(s^{-1}x)|\le1/L_r$
by \cref{lem:cutoff}\ref{it:cut-lip}, applied to the point $s^{-1}x$. Sum
over $r$. If all bits agree, the active prefixes, the index $n$ and the
selected field agree as well.
\end{proof}

\begin{figure}[ht]
\centering
\begin{tikzpicture}[x=0.62cm,y=0.75cm,font=\small,>=Latex]
 \fill[blue!8] (1.5,-0.35) rectangle (12.5,4.35);
 \fill[blue!25] (5.5,-0.35) rectangle (8.5,4.35);
 \node[blue!60!black,anchor=south] at (7,4.35) {$C_r$};
 \node[blue!50!black,anchor=south] at (2.6,4.35) {$E_r$};
 \draw[->] (-0.3,0)--(16,0);
 \draw[->] (0,-0.3)--(0,5) node[left] {$\phi_r$};
 \draw[dashed] (-0.1,4)--(15.6,4) node[right] {$1$};
 \draw[black!60] (0,0)--(2,0)--(3,1)--(4,2)--(5,3)--(6,4)--(8,4)--(9,3)--(10,2)--(11,1)--(12,0)--(15,0);
 \foreach \x/\y in {0/0,1/0,2/0,3/1,4/2,5/3,6/4,7/4,8/4,9/3,10/2,11/1,12/0,13/0,14/0,15/0}
   {\fill (\x,\y) circle (2pt);}
 \draw[red!70!black,thick] (-0.1,1.5)--(15.6,1.5) node[right] {$\omega$};
 \draw[red!70!black] (3,1) circle (4.5pt) (4,2) circle (4.5pt);
 \node[anchor=north] at (4,-0.1) {$x$};
 \draw[black!50] (3,-0.1)--(2.6,-0.6);
 \node[anchor=north east] at (2.75,-0.55) {$s^{-1}x$};
\end{tikzpicture}
\caption{The orbit cutoff $\phi_r$ along a path in the Schreier graph
(here $L_r=4$). It vanishes outside $E_r$ (light), changes by at most
$1/L_r$ along each edge, and equals $1$ on $C_r$ (dark). The bit $b_r$ is
one exactly where $\phi_r$ lies strictly above the threshold line at
height $\omega$. Two neighbours, such as $x$ and $s^{-1}x$ (circled), receive
different bits only if the line passes between their values. This
happens for $\omega$ in an interval of length at most $1/L_r$.}
\label{fig:threshold}
\end{figure}
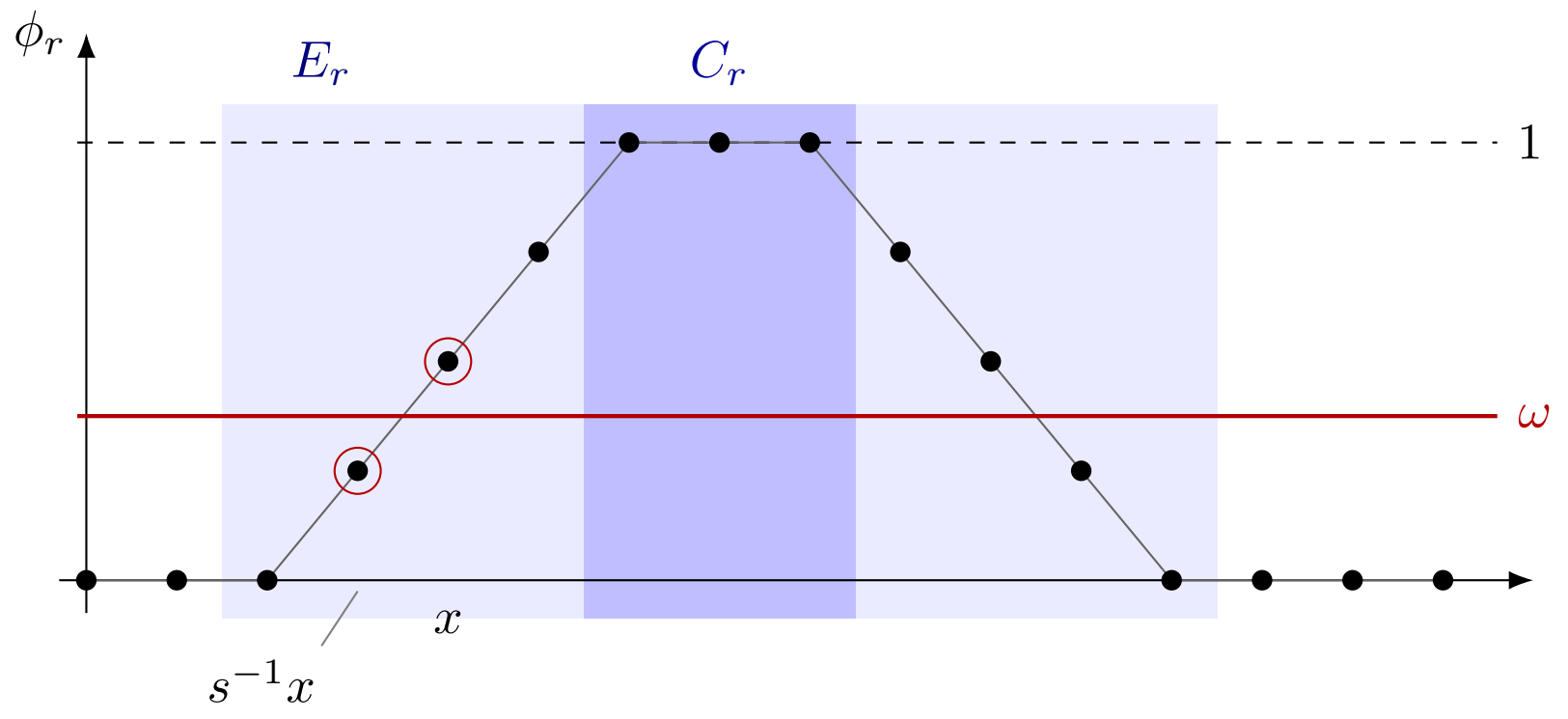

\Cref{fig:threshold} illustrates the common-threshold comparison.

\FloatBarrier

\subsection{The averaged field}\label{subsec:averaged}
At $x\in Y$ we use the selected field, evaluated at $x$, and average over
the threshold:
\begin{equation}\label{eq:F}
 F_s(x)=\begin{cases}\displaystyle\int_0^1c_s(x,\omega)(x)\,d\omega,&x\in Y,\\[1ex]
 0,&x\notin Y,\end{cases}\qquad(s\in S).
\end{equation}
The integrand is Borel on $X\times[0,1)$ by \cref{lem:field-measurable},
so each $F_s$ is a Borel function on $X$ with values in $[0,M]$.

\begin{proof}[Proof of \cref{thm:null}]
Let $F$ be the field defined by \eqref{eq:F}, and recall that
$M=4/\delta$. We claim that
\begin{equation}\label{eq:div-half}
 \Div F\ge\tfrac12\,\ind_Y\qquad\text{everywhere on }X.
\end{equation}
If $x\notin Y$, then also $s^{-1}x\notin Y$ for all $s$, because $Y$ is
invariant, so $\Div F(x)=0$. Now let $x\in Y$. Then every $s^{-1}x$ lies
in $Y$ as well. Suppressing $\omega$ from the notation, write
$c=c(x,\omega)$ and $c^{(s)}=c(s^{-1}x,\omega)$. Adding and subtracting
$c_s(s^{-1}x)$ inside the integral,
\begin{align*}
 \Div F(x)
 &=\int_0^1\sum_{s\in S}\Bigl(c_s(x)-c^{(s)}_s(s^{-1}x)\Bigr)\,d\omega\\
 &=\int_0^1\Div c(x)\,d\omega
 +\sum_{s\in S}\int_0^1\Bigl(c_s(s^{-1}x)-c^{(s)}_s(s^{-1}x)\Bigr)\,d\omega .
\end{align*}
All integrands are Borel in $\omega$ by \cref{lem:field-measurable}. By
\cref{lem:selected}\ref{it:sel-null}, $\Div c(x)>1$ for every $\omega$. In
the second term, if all bits of $x$ and $s^{-1}x$ agree, then
$c=c^{(s)}$ and the integrand vanishes. Otherwise the integrand is at
least $-M$, since all components take values in $[0,M)$. By
\cref{lem:selected}\ref{it:sel-couple} and~\eqref{eq:budget},
\[
 \Div F(x)\ge1-M\,|S|\sum_{r\ge1}\frac1{L_r}\ge\frac12 .
\]
This proves \eqref{eq:div-half}. Finally put $H=2F$. Then $H$ is a
nonnegative field supported on $Y$, with $\Div H\ge\ind_Y$ and $0\le
H_s\le2M=8/\delta$.
\end{proof}

The charge $\mu$ of \cref{thm:null} enters only through
\cref{lem:certificate}. All integrals above are ordinary Lebesgue
integrals over the threshold $\omega$.

\subsection{The converse}\label{subsec:converse}
The property that every $S$-invariant probability charge on $\Bor(X)$
vanishes on $m$-null Borel sets is in fact \emph{equivalent} to the
singular displacement gap, so \eqref{eq:H} is exactly the right
hypothesis.

\begin{corollary}\label{cor:converse}
Let $X$ be a $\sigma$-compact metrizable space, $S$ a finite nonempty
symmetric set of homeomorphisms of $X$, and $m$ a quasi-invariant Borel
measure on $X$ which is finite on compact sets. The following are equivalent.
\begin{enumerate}[label=\textup{(\roman*)}]
\item $(X,S,m)$ satisfies \eqref{eq:H} for some $\delta>0$.
\item Every $S$-invariant probability charge on $\Bor(X)$ vanishes on
all $m$-null Borel sets.
\end{enumerate}
\end{corollary}

\begin{proof}
(i)$\Rightarrow$(ii) is \cref{thm:null}. Suppose (i) fails. Then for
each $j\ge1$ there is a Borel probability measure $\nu_j\perp m$ with
$\max_s\|s_*\nu_j-\nu_j\|_\TV<1/j$. Let $N_j$ be an $m$-null Borel set
supporting $\nu_j$, and put
\[
 N=\bigcup_{j\ge1}N_j .
\]
This is an $m$-null Borel set with $\nu_j(N)=1$ for all $j$. Integration
against $\nu_j$ is an element of the closed unit ball of the dual of the
Banach space $\Bb(X)$ of bounded Borel functions with the supremum norm:
\[
 \nu_j(f)=\int f\,d\nu_j\qquad\bigl(f\in\Bb(X)\bigr).
\] By the Banach--Alaoglu theorem,
this ball is compact in the weak$^*$ topology, that is, the topology of
pointwise convergence on $\Bb(X)$. Hence the sequence $(\nu_j)$ has a
cluster point $\phi$ in this ball, which is the weak$^*$ limit of a
subnet $(\nu_{j_\alpha})$, with $j_\alpha\to\infty$. Since evaluation at a fixed function is
weak$^*$ continuous, $\phi$ is positive, $\phi(\ind_X)=1$ and
$\phi(\ind_N)=1$. Moreover, for every $f\in\Bb(X)$ with
$\|f\|_\infty\le1$ and every $s\in S$,
\[
 |\phi(f\circ s)-\phi(f)|=\lim_\alpha|\nu_{j_\alpha}(f\circ s)-\nu_{j_\alpha}(f)|
 \le\limsup_\alpha\|s_*\nu_{j_\alpha}-\nu_{j_\alpha}\|_\TV=0.
\]
Hence $\mu(A)=\phi(\ind_A)$, for $A\in\Bor(X)$, defines an $S$-invariant
probability charge with $\mu(N)=1$, and (ii) fails.
\end{proof}

\section{Isometric actions and the spherical theorem}\label{sec:sphere}

In this section we prove that a spectral gap forces uniqueness of the
invariant probability charge for isometric actions
(\cref{thm:gap-uniqueness}), and deduce \cref{thm:sphere}. There are two
steps. First we verify the singular displacement gap \eqref{eq:H}, so
that \cref{thm:null} applies and every invariant charge vanishes on null
sets. Smoothing a probability measure $\nu$ by an equivariant kernel of
small radius $r$ gives a density whose displacements are controlled by
those of $\nu$. The spectral gap, applied to the square root of this
density, shows that the density is close to the constant $1$ in $L^1$,
with a bound independent of $r$. Letting $r\to0$, the same bound holds
for $\|\nu-m\|_\TV$, which equals $2$ when $\nu$ is singular. Second, an invariant probability charge which vanishes on null sets
defines an invariant bounded functional $\Phi$ on $L^\infty(m)$, and a
duality argument shows that $\Phi(f)=\Phi(\ind_X)\,m(f)$, so the charge
is integration against $m$. For nonatomic $m$ this
second step is the classical uniqueness of invariant means on
$L^\infty(m)$ under a spectral gap~\cite[Theorem~1.1(a)]{LiNg}. We prove
it directly, because the duality argument is needed again, on Euclidean
balls, in \cref{sec:euclid}. It is therefore formulated abstractly.

\subsection{From spectral gap to the singular displacement gap}
\begin{proposition}\label{prop:isometric}
Let $X$ be a compact metric space and $m$ a Borel probability measure on
$X$ with full support. Let $S$ be a finite nonempty symmetric set of
isometries of $X$ preserving $m$, and suppose that $(X,S,m)$ has a
spectral gap with constant $\kappa$. Then every Borel probability measure
$\nu$ on $X$ satisfies
\[
 \|\nu-m\|_\TV\le\frac2\kappa\Bigl(\max_{s\in S}\|s_*\nu-\nu\|_\TV\Bigr)^{1/2}.
\]
In particular, $(X,S,m)$ satisfies \eqref{eq:H} with $\delta=\kappa^2$.
\end{proposition}

The proof uses an equivariant smoothing operator. Let $\varrho$ denote
the metric of $X$.

\begin{definition}[Smoothing]\label{def:smoothing}
In the setting of \cref{prop:isometric}, for $r>0$ define the
\emph{tent kernel} and its normalization
\[
 a_r(x,y)=\max\Bigl(0,1-\frac{\varrho(x,y)}r\Bigr),\qquad
 p_r(x,y)=\frac{a_r(x,y)}{\int_Xa_r(x,z)\,dm(z)}\qquad(x,y\in X).
\]
For a finite signed Borel measure $\eta$ on $X$, its \emph{smoothing}
$P_r\eta$ is the absolutely continuous measure
\[
 P_r\eta=h\,m,\qquad h(y)=\int_Xp_r(x,y)\,d\eta(x).
\]
\end{definition}

The denominator of $p_r$ is a continuous function of $x$, since $a_r$ is
$r^{-1}$-Lipschitz in $x$. It is positive, because $a_r(x,\cdot)\ge\frac12$
on the ball of radius $r/2$ about $x$, which has positive measure by the
assumption of full support. By compactness it is bounded below by a
positive constant. Thus $p_r$ is continuous and bounded on $X\times X$,
$\int p_r(x,y)\,dm(y)=1$ for every $x$, and $p_r(x,y)=0$ when
$\varrho(x,y)\ge r$. The following properties are what the proof needs.
Smoothing maps probability measures to absolutely continuous probability
measures, does not increase total variation, and commutes with the
generators (so
displacements of $\nu$ control displacements of the density), and it
recovers $\nu$ in the limit $r\to0$.

\begin{lemma}\label{lem:smoothing}
For every $r>0$, finite signed Borel measure $\eta$, and $s\in S$:
\begin{enumerate}[label=\textup{(\alph*)}]
\item\label{it:sm-prob} If $\eta$ is a probability measure, then so is $P_r\eta$.
\item\label{it:sm-contract} $|P_r\eta|(X)\le|\eta|(X)$.
\item\label{it:sm-equiv} $P_r(s_*\eta)=s_*(P_r\eta)$.
\item\label{it:sm-limit} If $\eta$ is a probability measure, then $\int
g\,d(P_r\eta)\to\int g\,d\eta$ as $r\to0$, for every $g\in C(X)$.
\end{enumerate}
\end{lemma}

\begin{proof}
\ref{it:sm-prob} follows from Tonelli's theorem and $\int
p_r(x,y)\,dm(y)=1$. For \ref{it:sm-contract}, $|h(y)|\le\int
p_r(x,y)\,d|\eta|(x)$, and Tonelli's theorem applied to the positive
measure $|\eta|$ gives $\int|h|\,dm\le|\eta|(X)$. For
\ref{it:sm-equiv}, since $s$ is an isometry we have
$a_r(sx,sy)=a_r(x,y)$, and since $s$ preserves $m$,
\[
 \int a_r(sx,z)\,dm(z)=\int a_r(sx,sz)\,dm(z)=\int a_r(x,z)\,dm(z).
\]
Hence $p_r(sx,sy)=p_r(x,y)$. The density of $P_r(s_*\eta)$ is therefore
\[
 \int p_r(sx,y)\,d\eta(x)=\int p_r(x,s^{-1}y)\,d\eta(x)=h(s^{-1}y),
\]
which is also the density of $s_*(hm)$, again because $s$ preserves $m$.
For \ref{it:sm-limit}, Fubini's theorem (the integrand is bounded) gives
$\int g\,d(P_r\eta)=\int g_r\,d\eta$ with $g_r(x)=\int g(y)p_r(x,y)\,dm(y)$.
Moreover $|g_r(x)-g(x)|\le\sup\{|g(y)-g(x)|:\varrho(x,y)<r\}$, which
tends to $0$ uniformly in $x$ by uniform continuity of $g$.
\end{proof}

\begin{proof}[Proof of \cref{prop:isometric}]
Put $D=\max_s\|s_*\nu-\nu\|_\TV$, fix $r>0$, and let $h$ be the density of
$P_r\nu$. By \cref{lem:smoothing}, $h\ge0$ and $\int h\,dm=1$. For each
$s\in S$ the density of $s_*(P_r\nu)=P_r(s_*\nu)$ is $h\circ s^{-1}$, so
\[
 \|h\circ s^{-1}-h\|_{L^1(m)}=\bigl|P_r(s_*\nu-\nu)\bigr|(X)\le\|s_*\nu-\nu\|_\TV\le D,
\]
using \cref{lem:TV-continuous}. Let $u=\sqrt h$, so that $\|u\|_{L^2}=1$.
The elementary inequality $(\sqrt a-\sqrt b)^2\le|a-b|$ $(a,b\ge0)$ gives
$\|u\circ s^{-1}-u\|_{L^2}^2\le\|h\circ s^{-1}-h\|_{L^1}\le D$. By the
spectral gap,
\[
 \kappa^2\bigl(1-m(u)^2\bigr)=\kappa^2\|u-m(u)\|_{L^2}^2\le D .
\]
By the Cauchy--Schwarz inequality,
\[
 \|h-1\|_{L^1}=\int|u-1|\,|u+1|\,dm\le\|u-1\|_{L^2}\|u+1\|_{L^2}
 =2\sqrt{1-m(u)^2}\le\frac{2\sqrt D}\kappa .
\]
Hence $|\int g\,d(P_r\nu)-\int g\,dm|\le2\sqrt D/\kappa$ for every
$g\in C(X)$ with $\|g\|_\infty\le1$. Letting $r\to0$ and using
\cref{lem:smoothing}\ref{it:sm-limit}, the same bound holds for
$|\int g\,d\nu-\int g\,dm|$, and \cref{lem:TV-continuous} gives the
inequality. If $\nu\perp m$ then $\|\nu-m\|_\TV=2$, so $D\ge\kappa^2$.
\end{proof}

\subsection{Invariant functionals and averaging}
The second step is a duality argument. We first convert the $L^2$
spectral gap into an inequality in $L^1$, by applying it to the
indicators of the superlevel sets of a function.

\begin{lemma}\label{lem:L1-gap}
Let $(\Omega,\lambda)$ be a probability space, and let $S$ be a finite
nonempty symmetric set of measure-preserving automorphisms of
$(\Omega,\lambda)$, such that
\[
 \max_{s\in S}\|u\circ s^{-1}-u\|_{L^2(\lambda)}\ge\kappa\,\|u-\lambda(u)\|_{L^2(\lambda)}
 \qquad\bigl(u\in L^2(\lambda)\bigr)
\]
for some $\kappa>0$. Then every $h\in L^1(\lambda)$ satisfies
\[
 \|h-\lambda(h)\|_{L^1}\le2\kappa^{-2}\sum_{s\in S}\|h-h\circ s\|_{L^1}.
\]
\end{lemma}

\begin{proof}
For $t\in\R$ consider the superlevel set and its measure,
\[
 A_t=\{h>t\},\qquad a_t=\lambda(A_t).
\]
Since $|h(x)-\lambda(h)|=|\int(h(x)-h(y))\,d\lambda(y)|\le\int|h(x)-h(y)|\,d\lambda(y)$,
the layer-cake formula $|h(x)-h(y)|=\int_\R|\ind_{A_t}(x)-\ind_{A_t}(y)|\,dt$
gives
\[
 \|h-\lambda(h)\|_{L^1}\le\iint|h(x)-h(y)|\,d\lambda(x)\,d\lambda(y)
 =\int_\R2a_t(1-a_t)\,dt .
\]
Apply the hypothesis to $u=\ind_{A_t}$. We have
$\|u-\lambda(u)\|_{L^2}^2=a_t(1-a_t)$, and, since $s$ preserves $\lambda$,
$\|\ind_{A_t}\circ s^{-1}-\ind_{A_t}\|_{L^2}^2=\lambda(sA_t\triangle A_t)=\lambda(A_t\triangle s^{-1}A_t)$.
Bounding the maximum by the sum,
\[
 \kappa^2a_t(1-a_t)\le\sum_{s\in S}\lambda(A_t\triangle s^{-1}A_t).
\]
Integrating over $t$ and using the layer-cake formula again in the form
$\int_\R\lambda(A_t\triangle s^{-1}A_t)\,dt=\|h-h\circ s\|_{L^1}$ gives
the claim. All interchanges of integrals have nonnegative integrands,
so Tonelli's theorem applies.
\end{proof}

By duality, an inequality of this kind forces invariant functionals on
$L^\infty$ to be invariant under averaging. We formulate this for a
bounded operator $E$ on $L^1(\lambda)$, which will be the mean
$h\mapsto\lambda(h)$ in this section and a spherical average in
\cref{sec:euclid}. For a $\sigma$-finite measure $\lambda$ we identify
$L^1(\lambda)^*$ with $L^\infty(\lambda)$,
and write $E^\dagger$ for the adjoint operator on $L^\infty(\lambda)$:
\[
 \int(E^\dagger f)\,h\,d\lambda=\int f\,(Eh)\,d\lambda\qquad\bigl(f\in L^\infty(\lambda),\ h\in L^1(\lambda)\bigr).
\]
\begin{definition}[Invariant functionals]\label{def:invariant-functional}
Let $S$ be a set of measure-preserving automorphisms of
$(\Omega,\lambda)$, so that $f\mapsto f\circ s$ is well defined on
$L^\infty(\lambda)$. A linear functional $\Phi$ on $L^\infty(\lambda)$ is
\emph{$S$-invariant} if
\[
 \Phi(f\circ s)=\Phi(f)\qquad\bigl(f\in L^\infty(\lambda),\ s\in S\bigr).
\]
\end{definition}

\begin{lemma}[Averaging invariant functionals]\label{lem:functional}
Let $(\Omega,\lambda)$ be a $\sigma$-finite measure space, and let $S$
be a finite nonempty symmetric set of measure-preserving automorphisms of
$(\Omega,\lambda)$. Let $E$ be a bounded linear operator on
$L^1(\lambda)$, and let $C\ge0$ with
\[
 \|h-Eh\|_{L^1}\le C\sum_{s\in S}\|h-h\circ s\|_{L^1}\qquad\bigl(h\in L^1(\lambda)\bigr).
\]
Then every bounded $S$-invariant linear functional $\Phi$ on
$L^\infty(\lambda)$ satisfies
\[
 \Phi(f)=\Phi(E^\dagger f)\qquad\bigl(f\in L^\infty(\lambda)\bigr).
\]
\end{lemma}

\begin{proof}
Give $(L^1(\lambda))^S$ the sum norm, and consider the difference
operator $T$:
\[
 \|(v_s)_{s\in S}\|=\sum_{s\in S}\|v_s\|_{L^1},\qquad
 T\colon L^1(\lambda)\to(L^1(\lambda))^S,\qquad
 Th=(h-h\circ s)_{s\in S}.
\]
Fix $f\in L^\infty(\lambda)$. On the subspace $T(L^1(\lambda))$
define
\[
 \Lambda(Th)=\int f\,(h-Eh)\,d\lambda .
\]
This is well defined. If $Th=Th'$, the hypothesis applied to
$h-h'$ gives $(h-h')-E(h-h')=0$. By the hypothesis again,
\[
 |\Lambda(Th)|\le\|f\|_{L^\infty}\|h-Eh\|_{L^1}\le C\|f\|_{L^\infty}\|Th\| .
\]
By the Hahn--Banach theorem, applied with the
sublinear functional $C\|f\|_{L^\infty}\|\cdot\|$, $\Lambda$ extends to a
bounded linear functional on $(L^1(\lambda))^S$. By the duality
$L^1(\lambda)^*=L^\infty(\lambda)$ in each coordinate, the extension has
the form $(v_s)\mapsto\sum_s\int G_sv_s\,d\lambda$ with $G_s\in
L^\infty(\lambda)$. For $h\in L^1(\lambda)$, invariance of $\lambda$
gives $\int G_s\,(h\circ s)\,d\lambda=\int(G_s\circ s^{-1})\,h\,d\lambda$,
so
\[
 \int\sum_{s\in S}\bigl(G_s-G_s\circ s^{-1}\bigr)\,h\,d\lambda=\Lambda(Th)
 =\int\bigl(f-E^\dagger f\bigr)\,h\,d\lambda .
\]
Since this holds for all $h\in L^1(\lambda)$, we get an equality of
classes
\[
 f-E^\dagger f=\sum_{s\in S}\bigl(G_s-G_s\circ s^{-1}\bigr)\qquad\text{in }L^\infty(\lambda).
\]
As
$\Phi$ is $S$-invariant and $S$ is symmetric,
$\Phi(G_s\circ s^{-1})=\Phi(G_s)$. Hence $\Phi(f)-\Phi(E^\dagger f)=0$.
\end{proof}

Let $\lambda$ be a measure and
$\mu$ a finite charge on the same measurable space $(\Omega,\mathcal A)$,
such that $\mu$ vanishes on all $\lambda$-null sets in $\mathcal A$. Any
two bounded $\mathcal A$-measurable representatives of a class in
$L^\infty(\lambda)$ differ only on a $\lambda$-null set, so by
\eqref{eq:integral-bound} they have the same integral against $\mu$.
We obtain a bounded linear
functional
\begin{equation}\label{eq:descend}
 \Phi_\mu\colon L^\infty(\lambda)\to\R,\qquad\Phi_\mu([f])=\mu(f),\qquad
 |\Phi_\mu([f])|\le\mu(\Omega)\|f\|_{L^\infty}.
\end{equation}
If $\mu$ is invariant under a measure-preserving automorphism $s$ of
$(\Omega,\lambda)$, then $\Phi_\mu$ is $s$-invariant by~\eqref{eq:integral-invariance}.

\subsection{Uniqueness under a spectral gap}
\begin{theorem}\label{thm:gap-uniqueness}
Let $X$ be a compact metric space and $m$ a Borel probability measure on
$X$ with full support. Let $S$ be a finite nonempty symmetric set of
isometries of $X$ preserving $m$, such that $(X,S,m)$ has a spectral
gap. Then $m$ is the only $S$-invariant probability charge on $\Bor(X)$.
\end{theorem}

\begin{proof}
Let $\kappa$ be a spectral gap constant and $\mu$ an $S$-invariant
probability charge on $\Bor(X)$. By \cref{prop:isometric}, $(X,S,m)$
satisfies \eqref{eq:H} with $\delta=\kappa^2$. Since $m$ is invariant,
it is quasi-invariant, so \cref{thm:null} applies: $\mu$ vanishes on
$m$-null Borel sets. Let $\Phi_\mu$ be the functional \eqref{eq:descend}
on $L^\infty(m)$. Apply \cref{lem:functional} with $(\Omega,\lambda)=(X,m)$,
$C=2\kappa^{-2}$, and the mean operator
\[
 E\colon L^1(m)\to L^1(m),\qquad Eh=m(h)\,\ind_X .
\]
Its hypothesis is \cref{lem:L1-gap}, and $E^\dagger f=m(f)\ind_X$. For a
Borel set $A$ we obtain
\[
 \mu(A)=\Phi_\mu([\ind_A])=\Phi_\mu\bigl(m(A)\ind_X\bigr)=m(A)\,\mu(X)=m(A).\qedhere
\]
\end{proof}

\subsection{Spectral gaps on spheres}
For $n\ge1$ let $\sigma_n$ be the normalized surface measure on the unit
sphere, which can be defined through Lebesgue measure $\Leb$ on
$\R^{n+1}$ by cones:
\[
 \mathbb S^n=\{x\in\R^{n+1}:|x|=1\},\qquad
 \sigma_n(A)=\frac{\Leb\{tx:x\in A,\ 0<t\le1\}}{\Leb\{y\in\R^{n+1}:0<|y|\le1\}}\qquad\bigl(A\in\Bor(\mathbb S^n)\bigr).
\]
It is the unique rotation-invariant Borel probability measure on
$\mathbb S^n$. Indeed, let $m_G$ be the Haar probability measure on
$G=\SO(n+1)$. For $f\in C(\mathbb
S^n)$ the function $x\mapsto\int f(gx)\,dm_G(g)$ is constant, by
transitivity and invariance of $m_G$. Integrating it against any
rotation-invariant probability measure $\nu$ and using Fubini's theorem
gives $\nu(f)=\int f(ge_1)\,dm_G(g)$, where $e_1$ is the first coordinate
vector, so $\nu$ is determined.

\begin{theorem}[Margulis, Sullivan, Drinfeld]\label{thm:gap-sets}
For every $n\ge2$ there is a finite nonempty symmetric set
$S\subseteq\SO(n+1)$ such that $(\mathbb S^n,S,\sigma_n)$ has a spectral
gap.
\end{theorem}

The existence of such sets is due to Margulis~\cite{Margulis}
and Sullivan~\cite{Sullivan} for $n\ge4$, and
to Drinfeld~\cite{Drinfeld} for $n=2,3$.

\begin{proof}[Proof of \cref{thm:sphere}]
Let $n\ge2$, and let $S\subseteq\SO(n+1)$ be a finite nonempty symmetric
set such that $(\mathbb S^n,S,\sigma_n)$ has a spectral gap. Rotations
are isometries of $\mathbb S^n$ preserving $\sigma_n$, which has full
support. By \cref{thm:gap-uniqueness}, $\sigma_n$ is the only
$S$-invariant probability charge on $\Bor(\mathbb S^n)$. This is the
second assertion of \cref{thm:sphere}. The first follows by applying it
to the set $S$ of \cref{thm:gap-sets}.
\end{proof}

\begin{remark}\label{rem:many-gap-sets}
Finite sets with a spectral gap are in fact plentiful, so the finite-set
form of \cref{thm:sphere} applies to many explicit families of
rotations. Benoist and de Saxc\'e~\cite[Theorem~1.2]{BdS} proved the
following. Let $G$ be a connected compact Lie group with simple Lie
algebra and Haar probability measure $m_G$. Let $S\subseteq G$ be a
finite nonempty symmetric set which generates a dense subgroup and such
that, in some fixed basis of the Lie algebra of $G$, every matrix
$\mathrm{Ad}(s)$, $s\in S$, has algebraic entries, where $\mathrm{Ad}$ is
the adjoint representation. Then
there is $\eps>0$ with
\[
 \Bigl\langle\frac1{|S|}\sum_{s\in S}u\circ s^{-1},\,u\Bigr\rangle\le(1-\eps)\|u\|_{L^2}^2
 \qquad\bigl(u\in L^2(m_G),\ m_G(u)=0\bigr).
\]
(For the special unitary group $G=\mathrm{SU}(d)$, the version with algebraic entries in the
defining representation is due to Bourgain and Gamburd~\cite[Theorem~1]{BG}.) Since $S$ is symmetric and $m_G$ is invariant,
\[
 \sum_{s\in S}\|u\circ s^{-1}-u\|_{L^2}^2=2|S|\Bigl(\|u\|_{L^2}^2-\Bigl\langle\frac1{|S|}\sum_{s\in S}u\circ s^{-1},u\Bigr\rangle\Bigr),
\]
and this gives \eqref{eq:spectral-gap} on $G$ with $\kappa=\sqrt{2\eps}$.
Pulling functions on $\mathbb S^n$ back to $G=\SO(n+1)$ along the orbit
map $g\mapsto ge_1$, which carries Haar measure to $\sigma_n$ and
intertwines the actions, transfers the gap to $(\mathbb
S^n,S,\sigma_n)$ whenever $\SO(n+1)$ has simple Lie algebra, that is,
for $n\ge2$, $n\ne3$. For $\mathbb S^3$ one can use $G=\mathrm{SU}(2)$
acting on $\mathbb S^3\cong\mathrm{SU}(2)$ by left multiplication, which
is an action by rotations.
\end{remark}

\section{Euclidean space}\label{sec:euclid}

In this section $d\ge3$ and $\sigma=\sigma_{d-1}$ is the normalized
surface measure on $\mathbb S^{d-1}$. We write $\Leb$ for Lebesgue
measure on $\R^d$, and $B(x,R)$ and $\bar B(x,R)$
for the open and closed balls:
\[
 B(x,R)=\{y:|y-x|<R\},\qquad\bar B(x,R)=\{y:|y-x|\le R\}.
\]
We work with the ring of bounded Borel sets and
the space of bounded Borel functions with bounded support:
\begin{align*}
 \Bor_b(\R^d)&=\{A\in\Bor(\R^d):A\text{ is bounded}\},\\
 B_c(\R^d)&=\{f\in\Bb(\R^d):f\text{ vanishes off a bounded set}\}.
\end{align*}

\subsection{Statement and plan}
For $a\in\R^d$ and $r\in\SO(d)$, let $r_a$ be the rotation $r$ about the
centre $a$. For $S\subseteq\SO(d)$ let $S_a$ be the corresponding set of
rotations about $a$:
\[
 r_a(x)=a+r(x-a),\qquad S_a=\{r_a:r\in S\}.
\]
Unlike the charges of \cref{sec:prelim}, the charges considered here are
defined on a ring rather than a $\sigma$-algebra and may a priori take
the value $\infty$.

\begin{definition}\label{def:local-charge}
A \emph{charge on bounded Borel sets} is a function
$\mu\colon\Bor_b(\R^d)\to[0,\infty]$ with $\mu(\varnothing)=0$ and
\[
 \mu(A\cup B)=\mu(A)+\mu(B)\qquad(A,B\in\Bor_b(\R^d),\ A\cap B=\varnothing).
\]
It is \emph{invariant} under a bijection $g$ of $\R^d$ preserving
$\Bor_b(\R^d)$ if
\[
 \mu(gA)=\mu(A)\qquad\bigl(A\in\Bor_b(\R^d)\bigr).
\]
\end{definition}

The following is the precise form of \cref{thm:euclid}.

\begin{theorem}\label{thm:euclid-precise}
Let $d\ge3$, let $p\ne q$ be points of $\R^d$, and let
$S^{(1)},S^{(2)}\subseteq\SO(d)$ be finite nonempty symmetric sets such
that $(\mathbb S^{d-1},S^{(i)},\sigma_{d-1})$ has a spectral gap for
$i=1,2$. Let $\mu$ be a charge on bounded Borel subsets of $\R^d$,
invariant under every element of $S^{(1)}_p\cup S^{(2)}_q$, such that
$\mu([0,1]^d)=1$. Then $\mu(A)$ is the Lebesgue measure of $A$ for every
bounded Borel set $A\subseteq\R^d$.
\end{theorem}

The first assertion of \cref{thm:euclid} is the case
$S^{(1)}=S^{(2)}=S$, where $S$ is given by \cref{thm:gap-sets} with
$n=d-1$.

Throughout the section we fix $d$, $p$, $q$, $S^{(1)}$ and $S^{(2)}$ as
in \cref{thm:euclid-precise}, with spectral gap constants $\kappa_1$ and
$\kappa_2$, and put
\[
 \ell=|p-q|,\qquad S=S^{(1)}_p\cup S^{(2)}_q,\qquad\Gamma=\langle S\rangle .
\]
The set $S$ is finite and symmetric, and $\Gamma$ is a countable group of
orientation-preserving isometries of $\R^d$. The proof follows the
outline in the introduction: the singular displacement gap for
$(\R^d,S,\Leb)$ (\cref{prop:euclid-H}), growth of $\mu$
(\cref{lem:dense}), null annihilation (\cref{lem:euclid-null}),
averaging invariance and countable additivity
(\cref{lem:average-invariance,lem:dini}), and a final mollification.

\subsection{Spherical averages}
\begin{definition}[Spherical averages]\label{def:spherical-average}
For $a\in\R^d$ and a bounded Borel function $f$ on $\R^d$, the
\emph{spherical average} of $f$ about $a$ is the function whose value at
$x$ is the mean of $f$ over the sphere about $a$ through $x$:
\[
 E_af(x)=\int_{\mathbb S^{d-1}}f\bigl(a+|x-a|\,\xi\bigr)\,d\sigma(\xi)\qquad(x\in\R^d).
\]
\end{definition}

The following properties keep us inside $B_c(\R^d)$, show that averaging
preserves Lebesgue integrals, and make $E_a$ self-adjoint. We use polar
coordinates about $a$: with $\omega_d$ the
surface area of $\mathbb S^{d-1}$, every Borel function $h\ge0$ on $\R^d$
satisfies
\begin{equation}\label{eq:polar}
 \int_{\R^d}h(x)\,dx=\omega_d\int_0^\infty\rho^{d-1}\int_{\mathbb S^{d-1}}h(a+\rho\xi)\,d\sigma(\xi)\,d\rho .
\end{equation}
The same formula holds for integrable $h$.

\Needspace{8\baselineskip}
\begin{lemma}\label{lem:averages}
Let $a\in\R^d$.
\begin{enumerate}[label=\textup{(\alph*)}]
\item\label{it:av-borel} If $f$ is bounded and Borel, then $E_af$ is
bounded and Borel, $\|E_af\|_\infty\le\|f\|_\infty$, and if $f$ vanishes
outside the closed ball $\bar B(a,R)$, so does $E_af$.
\item\label{it:av-integral} For every Borel function $h\ge0$, the formula
of \cref{def:spherical-average} defines a Borel function $E_ah\colon
\R^d\to[0,\infty]$ with $\int E_ah\,dx=\int h\,dx$. For a Borel
representative $h$ of a class in $L^1(\R^d)$, put
\[
 E_ah(x)=\begin{cases}\int_{\mathbb S^{d-1}}h(a+|x-a|\xi)\,d\sigma(\xi),&E_a|h|(x)<\infty,\\
 0,&\text{otherwise.}\end{cases}
\]
This induces a positive linear contraction of $L^1(\R^d)$ which preserves
integrals, and the same holds for $L^1(\bar B(a,R))$.
\item\label{it:av-symmetric} For every bounded Borel $f$ and every $h\in
L^1(\R^d)$,
\[
 \int(E_af)\,h\,dx=\int f\,(E_ah)\,dx .
\]
\item\label{it:av-continuous} If $f$ is continuous with compact support,
so is $E_af$.
\end{enumerate}
\end{lemma}

These are standard consequences of polar coordinates, and we prove them
in \cref{sec:averages}.

A single spherical average of a Borel function need not be continuous.
Averaging about a second centre smooths in the remaining directions,
except at that centre. A third average removes this last singularity.

\begin{lemma}\label{lem:three-averages}
For every $f\in B_c(\R^d)$, the function $E_pE_qE_pf$ is continuous and
compactly supported, and $\int E_pE_qE_pf\,dx=\int f\,dx$.
\end{lemma}

The proof, given in \cref{sec:averages}, rests on the classical formula
for the distribution of one coordinate on the sphere. On the sphere of
radius $t$ about $q$, the squared distance to $p$ is an affine function
of this coordinate, so after averaging about $p$ and then about $q$, the
result depends continuously on the distance to $q$ away from $q$ itself.
The third average removes this last singular point.

\subsection{The singular displacement gap}
We now show that the action of $\Gamma$ on $\R^d$, with Lebesgue measure
as reference measure, satisfies the singular displacement gap. The gap
is derived from the spectral gaps on spheres in three steps: an $L^1$
estimate for functions, its extension to measures by convolution, and a
smoothing argument with three averages.

The next lemma integrates the $L^1$ form of the spherical spectral gap
(\cref{lem:L1-gap}) over all spheres about a centre. It is used twice:
for the total variation estimate of \cref{prop:TV}, and for averaging on
balls in \cref{lem:average-invariance}.

\begin{lemma}\label{lem:polar-L1}
Let $a\in\R^d$, and let $Q\subseteq\SO(d)$ be a finite nonempty symmetric
set such that $(\mathbb S^{d-1},Q,\sigma)$ has a spectral gap with
constant $\kappa$. For every $R\in(0,\infty]$ and every $h\in
L^1(\bar B(a,R))$, where $\bar B(a,\infty)=\R^d$,
\[
 \|h-E_ah\|_{L^1}\le2\kappa^{-2}\sum_{r\in Q}\|h-h\circ r_a\|_{L^1}.
\]
\end{lemma}

\begin{proof}
Fix a Borel representative of $h$, and for $0<\rho<R$ consider its
restriction to the sphere of radius $\rho$,
\[
 h_\rho(\xi)=h(a+\rho\xi)\qquad(\xi\in\mathbb S^{d-1}).
\]
By \eqref{eq:polar}, $h_\rho\in L^1(\sigma)$ for almost every
$\rho\in(0,R)$. For such $\rho$, \cref{lem:L1-gap}, applied to $(\mathbb
S^{d-1},\sigma)$ and $Q$, gives
\[
 \int\bigl|h_\rho-\sigma(h_\rho)\bigr|\,d\sigma\le2\kappa^{-2}\sum_{r\in Q}\int|h_\rho-h_\rho\circ r|\,d\sigma .
\]
Here $\sigma(h_\rho)=E_ah(a+\rho\xi)$ for every $\xi$, and $h_\rho\circ
r(\xi)=(h\circ r_a)(a+\rho\xi)$. Multiply by $\omega_d\rho^{d-1}$ and
integrate over $\rho\in(0,R)$. By \eqref{eq:polar} the result is the
asserted inequality.
\end{proof}

The averaging operators act on finite signed Borel measures by duality:
the \emph{spherical average} $A_a\eta$ of $\eta$ is the finite signed
Borel measure
\begin{equation}\label{eq:measure-average}
 (A_a\eta)(B)=\int E_a\ind_B\,d\eta\qquad\bigl(B\in\Bor(\R^d)\bigr).
\end{equation}
Countable additivity follows from monotone convergence, applied to the
positive and negative parts of $\eta$. Since $E_a$ is positive with
$E_a1=1$, we have $\int f\,d(A_a\eta)=\int E_af\,d\eta$ for bounded Borel
$f$, and $|A_a\eta|(\R^d)\le|\eta|(\R^d)$. The next proposition
transfers \cref{lem:polar-L1} from densities to measures by
mollification.

\begin{proposition}\label{prop:TV}
Let $d\ge3$ and $a\in\R^d$, and let $Q\subseteq\SO(d)$ be a finite
nonempty symmetric set such that $(\mathbb S^{d-1},Q,\sigma_{d-1})$ has a
spectral gap with constant $\kappa$. Every Borel probability measure
$\nu$ on $\R^d$ satisfies
\[
 \|\nu-A_a\nu\|_\TV\le2\kappa^{-2}\sum_{r\in Q}\|\nu-(r_a)_*\nu\|_\TV .
\]
\end{proposition}

\begin{proof}
For $\eps>0$ let $\varphi_\eps=\ind_{B(0,\eps)}/\Leb(B(0,\eps))$, and let
$\nu_\eps=\nu*\varphi_\eps$ be the probability measure with density
\[
 h_\eps(x)=\int\varphi_\eps(x-y)\,d\nu(y).
\]
For every orientation-preserving isometry $g(x)=Ux+b$ ($U\in\SO(d)$) we
have $g_*(\nu*\varphi_\eps)=(g_*\nu)*\varphi_\eps$, because
$\varphi_\eps$ is invariant under $U$. Convolution with $\varphi_\eps$
does not increase total variation: for a finite signed Borel measure
$\eta$,
\[
 \int\Bigl|\int\varphi_\eps(x-y)\,d\eta(y)\Bigr|\,dx\le\int\!\!\int\varphi_\eps(x-y)\,dx\,d|\eta|(y)=|\eta|(\R^d).
\]
Hence, using $(r_a)^{-1}=(r^{-1})_a$ and the symmetry of $Q$,
\[
 \sum_{r\in Q}\|h_\eps-h_\eps\circ r_a\|_{L^1}
 =\sum_{r\in Q}\bigl|\nu_\eps-(r_a)_*\nu_\eps\bigr|(\R^d)
 \le\sum_{r\in Q}\|\nu-(r_a)_*\nu\|_\TV .
\]
By \cref{lem:averages}\ref{it:av-symmetric}, $A_a\nu_\eps$ has density
$E_ah_\eps$. So \cref{lem:polar-L1}, with $R=\infty$, gives
\[
 \bigl|\nu_\eps-A_a\nu_\eps\bigr|(\R^d)\le2\kappa^{-2}\sum_{r\in Q}\|\nu-(r_a)_*\nu\|_\TV .
\]
For every $\psi\in C_c(\R^d)$ we have
$\int\psi\,d\nu_\eps=\int\psi*\varphi_\eps\,d\nu$, where
$\psi*\varphi_\eps(y)=\int\psi(y+z)\varphi_\eps(z)\,dz$ (as $\varphi_\eps$
is even) converges to $\psi$ uniformly as $\eps\to0$. Let $f\in C_c(\R^d)$
with $\|f\|_\infty\le1$. Applying this to $\psi=f$ and to $\psi=E_af\in C_c(\R^d)$
(\cref{lem:averages}\ref{it:av-continuous}),
\[
 \Bigl|\int f\,d\nu-\int E_af\,d\nu\Bigr|
 =\lim_{\eps\to0}\Bigl|\int f\,d\nu_\eps-\int f\,d(A_a\nu_\eps)\Bigr|
 \le2\kappa^{-2}\sum_{r\in Q}\|\nu-(r_a)_*\nu\|_\TV .
\]
The left side is $|\int f\,d(\nu-A_a\nu)|$, and \cref{lem:TV-continuous}
completes the proof.
\end{proof}

The three-average lemma shows that averaging about $p$, $q$, $p$ turns
every measure into an absolutely continuous one.

\begin{lemma}\label{lem:abs-continuous}
For every Borel probability measure $\nu$ on $\R^d$, the measure
$A_pA_qA_p\nu$ is absolutely continuous with respect to Lebesgue
measure.
\end{lemma}

\begin{proof}
For bounded Borel $f$, $\int f\,d(A_pA_qA_p\nu)=\int E_pE_qE_pf\,d\nu$.
Let $N$ be a bounded Borel set with $\Leb(N)=0$. By
\cref{lem:three-averages}, $E_pE_qE_p\ind_N$ is continuous, nonnegative
and has integral $\Leb(N)=0$. Hence it vanishes identically, and
$(A_pA_qA_p\nu)(N)=0$. An arbitrary Lebesgue null Borel set is a
countable union of bounded ones.
\end{proof}

\begin{proposition}\label{prop:euclid-H}
Let $d\ge3$, let $p\ne q$ be points of $\R^d$, and let
$S^{(1)},S^{(2)}\subseteq\SO(d)$ be finite nonempty symmetric sets such
that $(\mathbb S^{d-1},S^{(i)},\sigma_{d-1})$ has a spectral gap with
constant $\kappa_i$ for $i=1,2$. Then Lebesgue measure $\Leb$ on $\R^d$ is
invariant under $S=S^{(1)}_p\cup S^{(2)}_q$, the set of rotations from
$S^{(1)}$ about $p$ and from $S^{(2)}$ about $q$, and
$(\R^d,S,\Leb)$ satisfies \eqref{eq:H} with constant
\[
 \delta=\frac2{2C_1+C_2},\qquad C_i=2|S^{(i)}|\kappa_i^{-2}.
\]
\end{proposition}

\begin{proof}
Rotations preserve Lebesgue measure. Let $\nu$ be a Borel probability
measure on $\R^d$ supported on a Lebesgue null set, and put
$D=\max_{s\in S}\|s_*\nu-\nu\|_\TV$. By \cref{prop:TV}, applied with
$(a,Q)=(p,S^{(1)})$ and $(a,Q)=(q,S^{(2)})$,
\[
 \|\nu-A_p\nu\|_\TV\le C_1D,\qquad\|\nu-A_q\nu\|_\TV\le C_2D .
\]
Since $A_p$ and $A_q$ do not increase total variation,
\[
 \nu-A_pA_qA_p\nu=(\nu-A_p\nu)+A_p(\nu-A_q\nu)+A_pA_q(\nu-A_p\nu)
\]
gives $\|\nu-A_pA_qA_p\nu\|_\TV\le(2C_1+C_2)D$. By
\cref{lem:abs-continuous}, the probability measures $\nu$ and
$A_pA_qA_p\nu$ are mutually singular, so their distance is two. Hence
$D\ge2/(2C_1+C_2)$.
\end{proof}

Lebesgue measure is invariant under $S$ and finite on compact sets.
Since $\R^d$ is $\sigma$-compact and metrizable, \cref{prop:euclid-H}
verifies the remaining hypothesis of \cref{thm:null}.

\subsection{Orbits and growth}
From now on $\mu$ is a charge as in \cref{thm:euclid-precise}. Since
invariance under the generators implies invariance under every
$g\in\Gamma$, we have $\mu(gA)=\mu(A)$ for all $g\in\Gamma$ and bounded
Borel $A$. The next lemma makes integration against $\mu$ possible on
bounded sets and supplies the growth bound needed to integrate
certificates. Its assertion about continuous invariant functions will
also identify the mollifications in the final proof.

\begin{lemma}\label{lem:dense}
Every continuous $\Gamma$-invariant function $h\colon\R^d\to\R$ is
constant. Consequently every $\Gamma$-orbit is dense in $\R^d$. Moreover
$\mu(A)<\infty$ for every bounded Borel set $A$, and there is a constant
$c_0$ such that
\begin{equation}\label{eq:growth}
 \mu\bigl(B(p,R)\bigr)\le c_0\,(R+1)^d\qquad(R\ge0).
\end{equation}
\end{lemma}

\begin{proof}
Let $h$ be continuous and $\Gamma$-invariant, and fix $\rho>0$. The
function $u(\xi)=h(p+\rho\xi)$ on $\mathbb S^{d-1}$ is continuous and
invariant under $S^{(1)}$. By the spectral gap, $u$ is
$\sigma$-almost everywhere equal to its mean, and since $\sigma$ has full
support and $u$ is continuous, $u$ is constant. So $h$ is constant on
every sphere about $p$, and likewise on every sphere about $q$. Hence
there is a function $w$ on $[0,\infty)$ with
\[
 h(x)=w(|x-p|)\qquad(x\in\R^d).
\]
Recall that $\ell=|p-q|$. For $n\ge1$, the distances to $p$ of the
points of the sphere of radius $n\ell$ about $q$ fill the interval
$[(n-1)\ell,(n+1)\ell]$, and $h$ is constant on that sphere. So $w$ is
constant on each of these intervals. They
overlap and cover $[0,\infty)$, so $w$, and hence $h$, is constant.

Let $x\in\R^d$, and let $F$ be the closure of its orbit $\Gamma x$. The
function $y\mapsto\dist(y,F)$ is continuous, and it is $\Gamma$-invariant
because $\Gamma$ acts by isometries preserving $F$. It vanishes at $x$,
hence everywhere. So $F=\R^d$.

Let $c=(\frac12,\dots,\frac12)$ be the centre of the unit cube. For
$y\in\R^d$ choose $g\in\Gamma$ with $|gc-y|<\frac14$. Then
\[
 B(y,\tfrac14)\subseteq B(gc,\tfrac12)=gB(c,\tfrac12)\subseteq g[0,1]^d,
 \qquad\text{so}\qquad\mu\bigl(B(y,\tfrac14)\bigr)\le\mu([0,1]^d)=1 .
\]
Every point of $\R^d$ lies within distance $\frac18$ of the lattice
$\Lambda=p+(4\sqrt d)^{-1}\Z^d$. Hence $B(p,R)$ is covered by the balls
of radius $\frac14$ about the points of $\Lambda\cap B(p,R+1)$, of which
there are at most $(8\sqrt d(R+1)+1)^d$. This gives \eqref{eq:growth}
with $c_0=(8\sqrt d+1)^d$, and finiteness on bounded sets follows.
\end{proof}

From now on $\mu$ is finite on bounded sets. For $f\in B_c(\R^d)$ we put
\[
 \mu(f)=\mu_B(f|_B),
\]
where $B$ is any bounded Borel set outside which $f$ vanishes and $\mu_B$
is the restriction of $\mu$ to Borel subsets of $B$, a charge in the
sense of \cref{def:charge}. This does not depend on $B$, since for
$B\subseteq B'$ the integrals against $\mu_B$ and $\mu_{B'}$ agree on
simple functions vanishing off $B$. The functional $f\mapsto\mu(f)$ is
linear and positive on $B_c(\R^d)$, satisfies the bounds
\eqref{eq:integral-bound} with $\Omega$ replaced by any bounded Borel set
outside which $f$ vanishes, and $\mu(f\circ g^{-1})=\mu(f)$ for every
$g\in\Gamma$.

\subsection{Null sets}
The charge $\mu$ is only defined on bounded sets, so a certificate for
an unbounded invariant Borel set $Y$ cannot be integrated over all of $\R^d$ as in
\cref{lem:certificate}. Instead we integrate it over the orbit
neighbourhoods $\Gamma_nK$ of a compact set $K$. The divergence
integrated over $\Gamma_nK$ is controlled by the charge of the shell
$\Gamma_nK\setminus\Gamma_{n-1}K$, which forces exponential growth
unless the charge of $Y\cap K$ vanishes.

\begin{lemma}\label{lem:growth-certificate}
Let $Y\subseteq\R^d$ be an invariant Borel set, let $M>0$, and let $H$
be a certificate for $Y$ with $0\le H_s\le M$ and $H_s=0$ off $Y$ for
all $s\in S$ (so $\Div H\ge\ind_Y$ everywhere). Let
$K\subseteq\R^d$ be compact, and suppose that the compact sets
$K_n=\Gamma_nK$ satisfy
\[
 \lim_{n\to\infty}\frac1n\log\bigl(1+\mu(K_n)\bigr)=0 .
\]
Then $\mu(Y\cap K)=0$.
\end{lemma}

\begin{proof}
Consider the charge of $Y$ inside the orbit neighbourhoods,
\[
 v_n=\mu(Y\cap K_n)\qquad(n\ge0).
\]
Let $n\ge1$ and $s\in S$. The functions
$\ind_{K_n}H_s$ and $\ind_{s^{-1}K_n}H_s$ lie in $B_c(\R^d)$, and
$\ind_{K_n}\,(H_s\circ s^{-1})=(\ind_{s^{-1}K_n}H_s)\circ s^{-1}$. By
invariance of $\mu$,
\[
 v_n\le\mu\bigl(\ind_{K_n}\Div H\bigr)
 =\sum_{s\in S}\mu\bigl(\ind_{K_n}H_s-\ind_{s^{-1}K_n}H_s\bigr)
 \le\sum_{s\in S}\mu\bigl(\ind_{K_n\setminus s^{-1}K_n}H_s\bigr).
\]
The last step drops a nonpositive term, using $H_s\ge0$. Since
$sK_{n-1}\subseteq K_n$, we have $K_{n-1}\subseteq s^{-1}K_n$. As
$0\le H_s\le M\ind_Y$,
\[
 v_n\le M|S|\,\mu\bigl(Y\cap(K_n\setminus K_{n-1})\bigr)=M|S|\,(v_n-v_{n-1}).
\]
With $\gamma=M|S|$ this rearranges to $\gamma\,v_{n-1}\le(\gamma-1)\,v_n$. Suppose
$v_0=\mu(Y\cap K)>0$. Taking $n=1$ shows $\gamma>1$, and then
$v_n\ge\theta^nv_0$ with $\theta=\gamma/(\gamma-1)>1$. Since $\mu(K_n)\ge v_n$,
this contradicts the hypothesis.
\end{proof}

\begin{lemma}\label{lem:euclid-null}
Every bounded Borel set $N\subseteq\R^d$ with $\Leb(N)=0$ satisfies
$\mu(N)=0$.
\end{lemma}

\begin{proof}
The saturation $Y=\Gamma N$ is an invariant Lebesgue null Borel set. By
\cref{prop:euclid-H,thm:null} it admits a nonnegative certificate
supported on $Y$ with components bounded by a constant $M$. Every
generator changes the distance to $p$ by at most $2\ell$: this is clear
for rotations about $p$, and for $s=r_q$ we have
$sx-p=r(x-p)+(I-r)(q-p)$ with $|(I-r)(q-p)|\le2\ell$. Thus
\[
 \bigl||sx-p|-|x-p|\bigr|\le2\ell\qquad(x\in\R^d,\ s\in S).
\]
Choose $R_0$ with $N\subseteq K=\bar B(p,R_0)$. Then $\Gamma_nK\subseteq\bar
B(p,R_0+2n\ell)$, and \eqref{eq:growth} gives
\[
 \mu(\Gamma_nK)\le c_0\,(R_0+2n\ell+2)^d ,
\]
a polynomial bound in $n$. By \cref{lem:growth-certificate}, $\mu(N)\le\mu(Y\cap
K)=0$.
\end{proof}

\subsection{Averaging and countable additivity}
Since $\mu$ vanishes on bounded null sets, it defines invariant bounded
functionals on $L^\infty$ of balls, to which the averaging lemma of
\cref{sec:sphere} applies. Recall that $B_c(\R^d)$ is the space of
bounded Borel functions with bounded support, and $E_a$ is the spherical
average about $a$.

\begin{lemma}\label{lem:average-invariance}
For $a\in\{p,q\}$ and $f\in B_c(\R^d)$ we have $\mu(f)=\mu(E_af)$.
\end{lemma}

\begin{proof}
Let $Q=S^{(1)}$ if $a=p$ and $Q=S^{(2)}$ if $a=q$, with spectral gap
constant $\kappa$. Choose $R>0$ such that $f$ vanishes outside
$B=\bar B(a,R)$. The restriction $\mu_B$ of $\mu$ to Borel subsets of $B$
is a finite charge, vanishing on Lebesgue null sets by
\cref{lem:euclid-null}, and invariant under the rotations $r_a$, $r\in
Q$, which preserve $B$ and Lebesgue measure. Let $\Phi$ be the bounded
$Q_a$-invariant functional on $L^\infty(B)$ given by \eqref{eq:descend}.
Apply \cref{lem:functional} with $\Omega=B$ and Lebesgue measure, the set
$Q_a$, the operator $E=E_a$ on $L^1(B)$, and $C=2\kappa^{-2}$. Its
hypothesis is \cref{lem:polar-L1}, and the adjoint $E^\dagger$ of $E$ for
the $L^1$--$L^\infty$ pairing is $E_a$ on $L^\infty(B)$ by
\cref{lem:averages}\ref{it:av-symmetric}. Since $E_af$ also vanishes
outside $B$ by \cref{lem:averages}\ref{it:av-borel}, we obtain
$\mu(f)=\Phi([f|_B])=\Phi([E_af|_B])=\mu(E_af)$.
\end{proof}

Three applications of \cref{lem:average-invariance} give
\begin{equation}\label{eq:L-average}
 \mu(f)=\mu(Pf),\qquad P=E_pE_qE_p\qquad\bigl(f\in B_c(\R^d)\bigr).
\end{equation}
Since $P$ maps indicators of bounded sets to continuous functions, this
upgrades $\mu$ to a measure.

\begin{lemma}\label{lem:dini}
For every bounded Borel set $B_0$, the restriction of $\mu$ to the Borel
subsets of $B_0$ is countably additive.
\end{lemma}

\begin{proof}
It suffices to show that $\mu(A_n)\to0$ for bounded Borel sets
$A_n\downarrow\varnothing$. The functions $f_n=P\ind_{A_n}$ are
continuous with compact support by \cref{lem:three-averages}, and
nonnegative. They decrease pointwise to zero, by dominated convergence
in each of the three spherical integrals. They all vanish outside the
compact set $K=\overline{\{f_1\ne0\}}$. They converge to zero uniformly
(Dini's theorem): given $\eps>0$, the sets $\{f_n\ge\eps\}$ are compact,
decrease, and have empty intersection, so $\{f_n\ge\eps\}=\varnothing$ for some $n$, and then
$f_k<\eps$ for all $k\ge n$. By \eqref{eq:L-average} and
\eqref{eq:integral-bound},
\[
 0\le\mu(A_n)=\mu(f_n)\le\|f_n\|_\infty\,\mu(K)\longrightarrow0 .\qedhere
\]
\end{proof}

\subsection{Proof of the Euclidean theorem}
\begin{proof}[Proof of \cref{thm:euclid-precise}]
By \cref{lem:dense,lem:dini}, $\mu$ is finite on bounded sets, and on
the Borel subsets of every bounded set it is a finite measure. We
mollify $\mu$. Fix a continuous radial probability density $\varphi_\eps$
on $\R^d$ supported on $B(0,\eps)$, $0<\eps\le1$, and put
\[
 h_\eps(y)=\int\varphi_\eps(y-x)\,d\mu(x)\qquad(y\in\R^d).
\]
For $y$ in a bounded set, the functions $x\mapsto\varphi_\eps(y-x)$ vanish
outside a common compact set $K$, and they depend uniformly continuously
on $y$. So
$|h_\eps(y)-h_\eps(y')|\le\mu(K)\sup_x|\varphi_\eps(y-x)-\varphi_\eps(y'-x)|$,
and $h_\eps$ is continuous. For $g\in S$, write $g(x)=Ux+b$ with
$U\in\SO(d)$. Then $\varphi_\eps(gy-gx)=\varphi_\eps(U(y-x))=\varphi_\eps(y-x)$,
so by invariance of $\mu$,
\[
 h_\eps(gy)=\int\varphi_\eps(gy-x)\,d\mu(x)=\int\varphi_\eps(gy-gx)\,d\mu(x)=h_\eps(y).
\]
By \cref{lem:dense}, $h_\eps$ is a constant $c_\eps$.

Let $J=\prod_{i=1}^d[\alpha_i,\beta_i]$ be a box, and let
$J'=\{x:\dist(x,J)\le1\}$, a bounded set on whose Borel subsets $\mu$ is
a finite measure. For $y\in J$ the function $x\mapsto\varphi_\eps(y-x)$
vanishes outside $J'$, and so does the integrated kernel
\[
 k_\eps(x)=\int_J\varphi_\eps(y-x)\,dy\qquad(x\in\R^d),
\]
which takes values in $[0,1]$. By Fubini's theorem for Lebesgue measure
on $J$ and the finite measure $\mu$ on $J'$,
\[
 c_\eps\Leb(J)=\int_Jh_\eps(y)\,dy=\int_{J'}k_\eps\,d\mu .
\]
As $\eps\to0$, $k_\eps(x)\to\ind_J(x)$ for every $x\notin\partial J$. The boundary
$\partial J$ is a bounded Lebesgue null set, so $\mu(\partial J)=0$ by
\cref{lem:euclid-null}. By dominated convergence, $c_\eps\Leb(J)\to\mu(J)$.
For $J=[0,1]^d$ this gives $c_\eps\to1$, and hence
\[
 \mu(J)=\Leb(J)\qquad\text{for every box }J .
\]
Boundaries of boxes are null for both measures, so $\mu$ and Lebesgue
measure also agree on half-open boxes. Fix a bounded half-open box $B_0$.
Finite disjoint unions of half-open boxes contained in $B_0$ form an
algebra generating the Borel subsets of $B_0$, and on it the two finite
measures agree. By the uniqueness of extensions, they agree on all Borel subsets of $B_0$.
Every bounded Borel set lies in such a box.
\end{proof}

\section{A characterization for isometric actions}\label{sec:complements}

For actions by isometries of a compact metric space preserving a Borel
probability measure of full support, the reference measure is the only invariant probability
charge on Borel sets exactly when there is a spectral gap. Sufficiency
is \cref{thm:gap-uniqueness}. Necessity rests on a classical theorem of
Schmidt about invariant means.

\begin{theorem}[Schmidt]\label{thm:schmidt}
Let $X$ be a compact metrizable space, $S$ a finite nonempty symmetric
set of homeomorphisms of $X$, and $m$ an $S$-invariant Borel probability
measure on $X$ without atoms such that the action is ergodic. If
integration against $m$ is the only $S$-invariant mean on $L^\infty(m)$,
that is, the only positive $S$-invariant linear functional $\phi$ on
$L^\infty(m)$ with $\phi(1)=1$, then $(X,S,m)$ has a spectral gap.
\end{theorem}

This is~\cite[Proposition~2.3]{Schmidt}, which is stated for a countable
group $\Gamma$ acting ergodically on a standard probability space without
atoms, here $\Gamma=\langle S\rangle$. Its conclusion is that there are a
finite set $Q\subseteq\Gamma$ and $\eps>0$ with
$\max_{g\in Q}\|u\circ g^{-1}-u\|_{L^2}\ge\eps\|u\|_{L^2}$ whenever
$m(u)=0$. If every $g\in Q$ is a word of length at most $L$ in $S$, the
triangle inequality and invariance of $m$ give
$\|u\circ g^{-1}-u\|_{L^2}\le L\max_{s\in S}\|u\circ s^{-1}-u\|_{L^2}$.
Applied to $u-m(u)$, this yields \eqref{eq:spectral-gap} with
$\kappa=\eps/L$. For modern treatments of this theorem
see~\cite[Section~5]{FS} and~\cite[Theorem~7.3]{BISG}.

\begin{theorem}\label{thm:isometric-criterion}
Let $X$ be a compact metric space, $m$ a Borel probability measure on $X$
with full support, and $S$ a finite nonempty symmetric set of isometries
of $X$ preserving $m$. Then $m$ is the only $S$-invariant probability
charge on $\Bor(X)$ if and only if $(X,S,m)$ has a spectral gap.
\end{theorem}

\begin{proof}
If $(X,S,m)$ has a spectral gap, $m$ is the only $S$-invariant
probability charge by \cref{thm:gap-uniqueness}.

Conversely, suppose that $m$ is the only $S$-invariant probability
charge on $\Bor(X)$, and let $\Gamma=\langle S\rangle$. We first show
that the action is ergodic and that integration is the only
$S$-invariant mean. If $A$ is a Borel set with $m(sA\triangle A)=0$ for all
$s$ and $0<m(A)<1$, then the conditional measure
\[
 \nu_A(B)=\frac{m(B\cap A)}{m(A)}\qquad(B\in\Bor(X))
\]
is an $S$-invariant probability measure with $\nu_A(A)=1\ne m(A)$, which
is impossible. If $\phi$ is an $S$-invariant mean on $L^\infty(m)$, then
\[
 \mu_\phi(A)=\phi([\ind_A])\qquad(A\in\Bor(X))
\]
is an $S$-invariant probability charge on $\Bor(X)$, hence equal to $m$.
So $\phi$ agrees with integration on simple functions. Simple functions are dense in $L^\infty(m)$ for the
essential supremum norm, and both functionals are continuous for it
(positivity applied to $\|f\|_{L^\infty}\pm f$ gives
$|\phi(f)|\le\|f\|_{L^\infty}$), so $\phi$ is integration.

If $m$ has no atoms, \cref{thm:schmidt} gives a spectral gap.

If $m$ has an atom $x_0$, its orbit consists of atoms of the same mass,
so it is finite. It is an invariant set of positive measure, so by
ergodicity it has measure one. Its complement is open and null, so full
support shows that $X=\Gamma x_0$ is finite and $\Gamma$ acts
transitively. If $X$ is a single point, \eqref{eq:spectral-gap} holds
trivially. Otherwise $L^2(m)$ is finite-dimensional, and
$u\mapsto\max_s\|u\circ s^{-1}-u\|_{L^2}$ is a continuous function on the
(nonempty, compact) unit sphere of the mean-zero subspace. It does not
vanish there, since an $S$-invariant function is constant by
transitivity. Its minimum is a spectral gap constant.
\end{proof}

Without full support, a spectral gap does not imply uniqueness. If
$X=\{a,b\}$, $S=\{\mathrm{id}\}$ and $m=\delta_a$, the spectral gap holds
trivially, but $\delta_b$ is invariant.

\appendix
\crefalias{section}{appendix}
\section{Spherical averages}\label{sec:averages}

In this appendix we prove the standard properties of the spherical
averages $E_a$ stated in \cref{sec:euclid}: \cref{lem:averages} and the
three-average smoothing \cref{lem:three-averages}. We keep the notation of
\cref{sec:euclid}. In particular $d\ge3$, $\sigma=\sigma_{d-1}$,
$\omega_d$ is the surface area of $\mathbb S^{d-1}$, and the polar
coordinate formula \eqref{eq:polar} holds.

\begin{proof}[Proof of \cref{lem:averages}]
\ref{it:av-borel} The map $(x,\xi)\mapsto f(a+|x-a|\xi)$ is Borel, so
its integral $E_af$ over $\xi$ is Borel. The rest is clear.

\ref{it:av-integral} Measurability follows as in \ref{it:av-borel}, using
monotone convergence for unbounded $h$. Since
$E_ah(a+\rho\xi)=\int h(a+\rho\eta)\,d\sigma(\eta)$ does not depend on
$\xi$, \eqref{eq:polar} gives the same value for $h$ and $E_ah$. If
$h$ is integrable, then $E_a|h|$ is integrable, so the set where it is
infinite is a Borel null set. If $h=h'$ almost everywhere, then
$E_a|h-h'|$ has integral zero, so $E_ah=E_ah'$ almost everywhere. Thus
$E_a$ acts on $L^1(\R^d)$, and $\int|E_ah|\le\int E_a|h|=\int|h|$. The
ball $\bar B(a,R)$ is a union of spheres about $a$, so $E_a$ also acts
on functions on it.

\ref{it:av-symmetric} By \eqref{eq:polar} and Fubini's theorem, both
sides equal
\[
 \omega_d\int_0^\infty\rho^{d-1}\Bigl(\int_{\mathbb S^{d-1}}f(a+\rho\eta)\,d\sigma(\eta)\Bigr)
 \Bigl(\int_{\mathbb S^{d-1}}h(a+\rho\xi)\,d\sigma(\xi)\Bigr)\,d\rho .
\]

\ref{it:av-continuous} If $x_j\to x$, then $f(a+|x_j-a|\xi)\to
f(a+|x-a|\xi)$ for every $\xi$, and dominated convergence gives
$E_af(x_j)\to E_af(x)$. Compact support follows from \ref{it:av-borel}.
\end{proof}

\begin{proof}[Proof of \cref{lem:three-averages}]
The last assertion follows from \cref{lem:averages}\ref{it:av-integral}
applied to $f^+$ and $f^-$. We use the classical formula for the
distribution of one coordinate on the sphere: for a unit vector $e$ and
bounded Borel $\phi$ on $[-1,1]$,
\begin{equation}\label{eq:latitude}
 \int_{\mathbb S^{d-1}}\phi(\langle\zeta,e\rangle)\,d\sigma(\zeta)=\int_{-1}^1\phi(u)\,g_d(u)\,du,
\end{equation}
where $g_d\colon\R\to[0,\infty)$ is the probability density
\[
 g_d(u)=c_d(1-u^2)^{(d-3)/2}\ \ (|u|<1),\qquad g_d(u)=0\ \ (|u|\ge1),\qquad
 c_d=\Bigl(\int_{-1}^1(1-u^2)^{(d-3)/2}\,du\Bigr)^{-1}.
\]
To verify \eqref{eq:latitude}, take $e$ to be the first basis vector and
compute $I=\int_{0<|x|\le1}\phi(x_1/|x|)\,dx$ in two ways, where
$\omega_j$ denotes the surface area of $\mathbb S^{j-1}\subseteq\R^j$. By
\eqref{eq:polar}, $I=\frac{\omega_d}d\int\phi(\zeta_1)\,d\sigma(\zeta)$.
By Fubini's theorem in the coordinates $x=(x_1,x')$, polar coordinates in
$x'\in\R^{d-1}$, and the substitution $x_1=\rho u$,
$|x'|=\rho\sqrt{1-u^2}$ (whose Jacobian is $\rho/\sqrt{1-u^2}$),
\[
 I=\omega_{d-1}\int_0^1\int_{-1}^1\phi(u)\,\rho^{d-2}(1-u^2)^{(d-2)/2}\,\frac{\rho}{\sqrt{1-u^2}}\,du\,d\rho
 =\frac{\omega_{d-1}}d\int_{-1}^1\phi(u)(1-u^2)^{(d-3)/2}\,du .
\]
Comparing the two expressions and taking $\phi=1$ to identify the
constant gives \eqref{eq:latitude}. Only $g_d\in L^1(\R)$ will be used.

The first average is a function of the distance to $p$:
\[
 E_pf(x)=k(|x-p|),\qquad k(\rho)=\int_{\mathbb S^{d-1}}f(p+\rho\xi)\,d\sigma(\xi)\quad(\rho\ge0),
\]
where $k$ is a bounded Borel function on $[0,\infty)$. Recall that
$\ell=|p-q|$. Let $e=(q-p)/\ell$, and let $H(v)=k(\sqrt v)$ for $v\ge0$ and $H(v)=0$ for
$v<0$, a bounded Borel function on $\R$. For a point $x=q+t\eta$ with
$t>0$ and $\eta\in\mathbb S^{d-1}$, the squared distance from $p$ to a
point of the sphere about $q$ through $x$ is
$|q+t\zeta-p|^2=a_t+b_t\langle\zeta,e\rangle$, where
\[
 a_t=\ell^2+t^2,\qquad b_t=2\ell t .
\]
So by \eqref{eq:latitude} and the substitution $v=a_t+b_tu$,
\[
 E_qE_pf(x)=\int_{-1}^1H(a_t+b_tu)\,g_d(u)\,du
 =\int_\R H(v)\,g_{(t)}(v)\,dv,\qquad g_{(t)}(v)=\frac1{b_t}\,g_d\Bigl(\frac{v-a_t}{b_t}\Bigr).
\]
Affine maps act continuously on $L^1(\R)$: if $(a,b)\to(a_0,b_0)$ with
$b_0>0$, then $b^{-1}g((\cdot-a)/b)\to b_0^{-1}g((\cdot-a_0)/b_0)$ in
$L^1(\R)$ for every $g\in L^1(\R)$. For continuous $g$ with compact
support this follows from uniform continuity and a common compact
support. The general case follows by approximation in $L^1$, since
continuous functions with compact support are dense in $L^1(\R)$
 and these affine changes preserve the
$L^1$ norm. Hence $t\mapsto g_{(t)}$ is continuous from $(0,\infty)$ into
$L^1(\R)$, and $E_qE_pf(q+t\eta)$, which depends only on $t$, is
continuous in $t>0$. Thus $F=E_qE_pf$ is bounded and continuous on
$\R^d\setminus\{q\}$.

Finally, let $x_j\to x$ in $\R^d$. For each $\xi\in\mathbb S^{d-1}$,
$p+|x_j-p|\xi\to p+|x-p|\xi$, and the limit equals $q$ for at most one
$\xi$ (for none if $x=p$, since $p\ne q$). Since $\sigma$ has no atoms
and $F$ is bounded and continuous off $q$, dominated convergence gives
$E_pF(x_j)\to E_pF(x)$. Hence $E_pE_qE_pf$ is continuous. By
\cref{lem:averages}\ref{it:av-borel} it vanishes outside a bounded set,
so it is compactly supported.
\end{proof}

\end{document}